\documentclass[11pt,reqno,oneside]{amsart}
\usepackage{amscd}
\usepackage{amsmath}
\usepackage{latexsym}
\usepackage{amsfonts}
\usepackage{amssymb}
\usepackage{amsthm}
\usepackage{graphicx}
\usepackage{verbatim}
\usepackage{mathrsfs}
\usepackage{enumerate}
\usepackage{hyperref}
\usepackage{fullpage}
\usepackage{xcolor}
\usepackage[latin1]{inputenc}

\theoremstyle{plain}
\theoremstyle{plain}\newtheorem{theorem}{Theorem}[section]
\theoremstyle{plain}\newtheorem{lemma}[theorem]{Lemma}
\theoremstyle{plain}\newtheorem{coro}[theorem]{Corollary}
\theoremstyle{plain}\newtheorem{proposition}[theorem]{Proposition}
\theoremstyle{plain}\newtheorem{remark}{Remark}[section]

\numberwithin{equation}{section}

\newcommand{\norm}[1]{\left\|#1\right\|}

\newcommand{\R}{\mathbb{R}}
\newcommand{\be}{\begin{equation}}
\newcommand{\ee}{\end{equation}}
 \newcommand{\ba}{\begin{aligned}}
 \newcommand{\ea}{\end{aligned}}

  \newcommand{\ben}{\begin{enumerate}}
   \newcommand{\een}{\end{enumerate}}

\newcommand{\Rmnum}[1]{\expandafter\@slowromancap\romannumeral #1@}

\allowdisplaybreaks

\makeatletter
\@namedef{subjclassname@2020}{\textup{2020} Mathematics Subject Classification}
\makeatletter

\begin{document}
\title{On a 1d nonlocal transport of the incompressible porous media equation}
\author{Caifeng Liu$^{1}$,\,Wanwan Zhang$^{2*}$}
\address{$^1$ Department of Mathematics, Jiangxi University of Finance and Economics, Nanchang 330013, Jiangxi, P. R. China}
\email{liucaifeng@jxufe.edu.cn}

\address{$^2$ School of Mathematics and Statistics, Jiangxi Normal University, Nanchang, 330022, Jiangxi, P. R. China}
\email{zhangww@jxnu.edu.cn}

\thanks{$^*$Corresponding author}

\begin{abstract}
Recently, Kiselev and Sarsam proposed the following nonlocal transport equation as a one-dimensional analogue of the 2D incompressible porous media (IPM) equation
\begin{eqnarray*}
\partial_t\rho+u\partial_x\rho= 0,~u=gH_a\rho,
\end{eqnarray*}
where the transform $H_a$ is defined by
\begin{eqnarray*}
H_af(x)=\frac{1}{\pi}P.V.\int\limits_{\mathbb{R}}\frac{a^2f(y)}{(x-y)((x-y)^2+a^2)}dy.
\end{eqnarray*}
In the work Kiselev-Sarsam (2025) \cite{[Kiselev-Sarsam]}, the authors proved the local well-posedness for this 1D periodic IPM model as well as finite time blow-up for a class of smooth initial data. In this paper, we present several new weighted inequalities for the transform $H_a$ in the setting of the real line. Based on these integral inequalities, we also prove the finite time blow-up for this 1D IPM model on the real line.
\end{abstract}
\smallskip
\maketitle
\section{Introduction}
In this paper, we consider the Cauchy problem of the following one-dimensional transport equation with a nonlocal velocity
\begin{equation}\label{KS}
\left\{\ba
&\partial_t\rho+u\partial_x\rho= 0, ~(x,t)\in \R\times\R_+,\\
&u=gH_a\rho,\\
&\rho(x,0)=\rho_{0}(x),~x\in\R, \ea\ \right.
\end{equation}
where $a,g>0$ are fixed constants, and the transform $H_a$ is defined by
\begin{eqnarray}\label{definition-Ha}
H_af(x)=\frac{1}{\pi}P.V.\int\limits_{\R}\frac{a^2f(y)}{(x-y)((x-y)^2+a^2)}dy.
\end{eqnarray}
The non-local active scalar transport equation \eqref{KS} was recently proposed by Kiselev and Sarsam in \cite{[Kiselev-Sarsam]} as a one-dimensional model analogy for two-dimensional incompressible porous media (IPM) equation, which models the transport of a scalar density $\rho=\rho(x,t)$ by an incompressible fluid velocity field $u(x,t)=(u_1(x,t),u_2(x,t))$ under the effects of Darcy's law and gravity.
As described in \cite{[Kiselev-Sarsam]}, the motivation for deriving the equation \eqref{KS} is to model the evolution of the boundary trace of a solution to the initial value problem with a no-flux boundary condition for the 2D IPM in the periodic half-plane $\Omega=\mathbb{T}\times(0,\infty)$ with $\mathbb{T}=[-\pi,\pi)$ denoting the circle:
\begin{equation}\label{2-d-IPM}
\left\{\ba
&\partial_t\rho+u\cdot\nabla\rho= 0, ~(x,t)\in \Omega\times\R_+,\\
&u=-\nabla P-(0,g\rho),~\nabla\cdot u=0,\\
&u_2=0,~(x,t)\in \partial\Omega\times\R_+,\\
&\rho(\cdot,0)=\rho_{0},~x\in\Omega, \ea\ \right.
\end{equation}
where $P=P(x,t)$ is the scalar pressure, $g>0$ is the constant of gravitational acceleration, and $\rho_{0}=\rho_{0}(x)$ is the initial density.
In \eqref{KS}, $g>0$ is then the constant of gravitational acceleration, $\rho=\rho(x,t)$ is the scalar boundary density,  $\rho_0=\rho_0(x)$ is the initial boundary density, and $u=u(x,t)$ is the scalar 1D fluid velocity. One can refer to \cite{[Kiselev-Sarsam]} for more details of the derivation of \eqref{KS} from the IPM equation \eqref{2-d-IPM}. We also refer the readers to \cite{[Bianchini-Cordoba-Martninez],[Cordoba-Martninez]} and references therein for recent results on the 2D IPM equation.
In \cite{[Kiselev-Sarsam]}, the authors also remarked in detail how the 1D IPM equation \eqref{KS} parallels the well-known C\'{o}rdoba-C\'{o}rdoba-Fontelos (CCF) equation:
\begin{eqnarray}\label{CCF}
\partial_t\theta-H\theta\partial_x\theta=0,
\end{eqnarray}
where $H$ is the Hilbert transform defined by
\begin{eqnarray*}
H\theta(x)=\frac{1}{\pi}P.V.\int\limits_{\R}\frac{\theta(y)}{x-y}dy.
\end{eqnarray*}
The CCF equation and its natural generalizations have been extensively studied.
In \cite{[Cordoba-Cordoba-Fontelos05]}, C\'{o}rdoba-C\'{o}rdoba-Fontelos first showed the finite time singularity formation of \eqref{CCF} for a class of smooth even initial data.
Specifically, the blow-up proof in \cite{[Cordoba-Cordoba-Fontelos05]} is based on an ingenious inequality:
\begin{eqnarray}\label{CCF-inequality}
-\int\limits^\infty_0\frac{H f(x) f'(x)}{x^{1+\sigma}}dx
\geq
C_{\sigma}\int\limits^\infty_0\frac{(f(x))^2}{x^{2+\sigma}}dx,
\end{eqnarray}
where $-1<\sigma<1$, $C_\sigma>0$ is a constant depending only on $\sigma$, and $f$ is an even bounded smooth function on $\R$. The proof of \eqref{CCF-inequality} in \cite{[Cordoba-Cordoba-Fontelos05]} uses the Meillin transform and complex analysis.
The inequality \eqref{CCF-inequality} was later used in \cite{[Li-Rodrigo08]} to show the finite time blow-up of smooth solutions to the following fractionally dissipative CCF equation:
\begin{eqnarray}\label{dissipative-CCF}
\partial_t\theta-H\theta\partial_x\theta+(-\Delta)^{\frac{\gamma}{2}}\theta=0,~0<\gamma<\frac12.
\end{eqnarray}
The well-posedness for \eqref{dissipative-CCF} with $0<\gamma<2$  in the critical Sobolev space was studied in \cite{[Dong2008]}.
Many different proofs of the finite time blow-up for \eqref{CCF} and \eqref{dissipative-CCF} were provided in \cite{[Kiselev],[Li-Rodrigo20],[Silvestre-Vicol]}.
In particular, Silvestre and Vicol in \cite{[Silvestre-Vicol]} provided four elegant blow-up proofs for the CCF equation \eqref{CCF}.
Based on completely real variable arguments,
Kiselev \cite{[Kiselev]} showed the following inequality: for any even bounded $C^1$ function $f$ with $f(0)=0$ and $f'\geq0$ for $x>0$,
\begin{eqnarray}\label{nonlinear-inequality-H-Kiselev}
-\int\limits^1_0\frac{H f(x) f'(x)(f(x))^{p-1}}{x^{\sigma}}dx
\geq
C_{\sigma,p}\int\limits^1_0\frac{(f(x))^{p+1}}{x^{1+\sigma}}dx,
\end{eqnarray}
where $p\geq1$ and $\sigma>0$.
In this direction, Li and Rodrigo in \cite{[Li-Rodrigo20]} gave several new, elementary and real variable based proofs of the inequalities \eqref{CCF-inequality} and \eqref{nonlinear-inequality-H-Kiselev}. The authors there also established the following integral inequality for the Hilbert transform $H$: for any even Schwartz function $f$,
\begin{eqnarray}\label{nonlinear-inequality-exponential-weight}
-\int\limits^\infty_0f'(x)Hf(x)\frac{e^{-x}}{x}dx
\geq \frac{1}{2\pi}
\int\limits\limits^\infty_0\frac{(f(0)-f(x))^2}{x^2}dx
-1000\|f\|^2_{L^\infty},
\end{eqnarray}
which can be applied to prove the finite time blow-up of smooth solutions to \eqref{CCF} and \eqref{dissipative-CCF}.
However, the question of finite time blow-up for \eqref{dissipative-CCF} with $\frac12\leq\gamma<1$ in the class of smooth solutions is still open.
It should be mentioned that Hoang-Radosz in \cite{[Hoang-Radosz]} introduced and studied a nonlocal active scalar, inspired by the CCF equation \eqref{CCF}, and proved that either a cusp- or needle- like singularity forms in finite time.
At last, the natural generalizations of the one-dimensional Hilbert-type models \eqref{CCF} and \eqref{dissipative-CCF} to multi-dimensions were the Riesz-type models, which also have been intensively studied in the literature (see, e.g.,\cite{[Balodis-Cordoba],[Dong],[Dong-Li],[Jiu-Zhang],[Li-Rodrigo09],[Li-Rodrigo20],[Zhang]}).

Finally, in \cite{[Kiselev-Sarsam]}, the authors also proved the local well-posedness for the equation \eqref{KS} posed on the circle and adapted the arguments applied to the Hilbert transform $H$ in \cite{[Kiselev]} to show that for any $a,\sigma>0$, it holds that
\begin{eqnarray}\label{IPM-inequality}
-\int\limits_0^{\frac\pi2}\frac{H_af(x)f'(x)}{x^\sigma}dx\geq C_{a,\sigma}\int\limits_0^{\frac\pi2}\frac{(f(x))^2}{x^{1+\sigma}}dx,
\end{eqnarray}
where $C_{a,\sigma}$ is a universal constant depending only on $\sigma$ and $a$, and $f$ is an even and nonnegative smooth function defined on $\mathbb{T}$  with $f(0)=0$ and $f'\geq0$ on $[0,\pi)$. Based on the integral inequality \eqref{IPM-inequality}, Kiselev and Sarsam proved the finite time singularity of solutions to \eqref{KS} for a class of  smooth even initial data in the setting of the periodic circle.

The purpose of this paper is to revisit the Kiselev-Sarsam equation \eqref{KS} and present several weighted integral inequalities for the transform $H_a$ in the setting of the real line. Based on these integral inequalities, we also prove the finite time blow-up for this 1D IPM model on the real line.

The remaining part of this paper is organized as follows. In Section 2, we first present the local well-posedness for the model \eqref{KS}, and then recall some properties of smooth solutions to \eqref{KS}.
In Section 3 we first derive a point-wise inequality (see Proposition \ref{lower-bound-Ha}) for the transform $H_a$ acting on even and non-increasing (on $(0,\infty)$) functions on $\R$, and then establish the C\'{o}rdoba-C\'{o}rdoba-Fontelos-type integral inequality for the transform $H_a$ (see Proposition \ref{nonlinear-inequality-Ha-prop-3.2}) by a simple application of Fubini's theorem and integration by parts. In this section we also prove the weighted integral inequality for the transform $H_a$ acting on even but non-monotone (on $(0,\infty)$) functions on $\R$ (see Proposition \ref{nonlinear-weighted-inequality-exponential-prop-3.6}), and then prove the finite time blow-up of smooth solution to \eqref{KS} for a class of smooth even initial data.
In Section 4 we give another proof of the finite time blow-up via the method of telescoping sums. This proof is based on a local in space upper bound for the nonlinearity (see Lemma \ref{upper-bound-nonlinear-vicol}). Finally, Section 5 is devoted to the proof of the Kiselev-type inequality for the transform $H_a$.

Throughout this paper, we will use $C$ to denote a positive constant, whose value may change from line to line, and write $C_{a,\sigma}$ to emphasize the dependence of a constant on $a$ and $\sigma$. The symbol $i$ denotes the  complex number such that $i^2=-1$.
For $p\in[1,\infty]$, we denote $L^p(\R)$ the standard $L^p$-space and its  norm by  $\|\cdot\|_{L^p(\R)}$.
For $s\geq0$, we use the notation $H^s(\R)$ to denote the nonhomogeneous Sobolev space of $s$ order,
whose endowed norm is denoted by $\|\cdot\|_{H^s(\R)}=\|\cdot\|_{L^2(\R)}+\|(-\Delta)^{\frac s2}(\cdot)\|_{L^2(\R)}$,
where the fractional Laplacian $(-\Delta)^{\frac s2}$ is defined through the Fourier transform
\begin{eqnarray*}
\widehat{(-\Delta)^{\frac s2}f}(\xi)=(2\pi|\xi|)^s\widehat{f}(\xi).
\end{eqnarray*}
For a sake of the convenience, the $L^{p}(\R)$-norm of a function $f$ is sometimes abbreviated as $\|f\|_{L^p}$ and
the $H^s(\R)$-norm as $\norm{f}_{H^s}$. Finally, $BMO(\mathbb{R})$ denotes the classical $BMO$ function space on $\mathbb{R}$ with the norm notation $\|\cdot\|_{BMO(\R)}$.
\section{Preliminaries}
In this section, we now present the local well-posedness for the model \eqref{KS} and recall some properties of the solutions to \eqref{KS}.

We begin with the local well-posedness of \eqref{KS} in the Sobolev space $H^s(\R)$ for $s>\frac32$, whose proof will be sketched in the appendix for completeness.
\begin{proposition}\label{local-well-posedness}
Fix $a,g>0$ and assume that $\rho_0\in H^s(\R)$ for $s>\frac32$. Then there exists $T>0$ such that the equation \eqref{KS} has a unique solution $\rho$ in $C([0,T); H^s(\R))\cap {Lip}((0,T);H^{s-1}(\R)).$ Furthermore if $T^\ast$ is the first blow-up time such that $\rho$ cannot be continued in $C([0,T^\ast); H^s(\R))$, then we must have that
\begin{eqnarray*}
\displaystyle\limsup_{t\rightarrow T^\ast}\|\rho(t)\|_{H^s(\R)}=\infty~{\rm if~and~only~if}~\lim_{t\rightarrow T^\ast}\int\limits_0^{t}\|\partial_x\rho(\cdot,s)\|_{L^\infty(\R)}ds=\infty.
\end{eqnarray*}
\end{proposition}
Next we recall that the even symmetry and the monotonicity  of initial data are preserved by the solution to \eqref{KS}, whose proof can be consulted in \cite{[Kiselev-Sarsam]}.
\begin{lemma}\label{nondecreasing-property}
Fix $a,g>0$ and assume $\rho_0$ is an smooth even function on $\R$ and non-increasing on $(0,\infty)$. Let $\rho$ be a smooth, decaying at infinity solution to the problem \eqref{KS}. Then $\rho(\cdot,t)$ is even and non-increasing on $(0,\infty)$.
\end{lemma}
If $f$ is an even function, the formula \eqref{definition-Ha} becomes
\begin{eqnarray}\label{velocity-representation}
\begin{split}
H_af(x)
&=\frac{2a^2x}{\pi}\int\limits^\infty_0\frac{(f(y)-f(x))(x^2+3y^2+a^2)}{(x^2-y^2)((x-y)^2+a^2)((x+y)^2+a^2)}dy.
\end{split}
\end{eqnarray}
Note that if $\rho$ is even on $\R$ and monotone non-increasing on $(0,\infty)$, then $H_a\rho(x)>0$ at every point $x>0$ unless $\rho$ is constant.

We proceed to state the $L^\infty$ maximum principle, which follows from the transport structure of \eqref{KS}. One can refer to \cite{[Kiselev-Sarsam]} for more details of the proof.
\begin{lemma}\label{maximum-principle}
Fix $a,g>0$ and suppose for $s>\frac72$ that $\rho$ is a solution in $H^s(\mathbb{R})$ to \eqref{KS} on the time interval $[0,T)$. Then, given $\rho_0\in H^s(\mathbb{R})\subseteq L^\infty(\mathbb{R})$, we have $\|\rho(t)\|_{L^\infty(\mathbb{R})}=\|\rho_0\|_{L^\infty(\mathbb{R})}$.
\end{lemma}

\section{Several weighted integral inequalities for the transform $H_a$}
\subsection{Monotone decreasing case}
We begin with a pointwise inequality for the transform $H_a$ acting on even and non-increasing (on $(0,\infty)$) functions on $\R$.
The similar inequality for the Hilbert transform $H$ was established in \cite{[Li-Rodrigo20]}.
\begin{proposition}\label{lower-bound-Ha}
Let $f:\R\rightarrow\R$ be an even Schwartz function
and monotone decreasing on $(0,\infty)$. Then for any $a>0$ and $x>0$, we have
\begin{eqnarray*}
H_a f(x)\geq\frac{2}{\pi}\cdot\frac{a^2}{x(x^2+a^2)}\int\limits^x_0(f(y)-f(x))dy .
\end{eqnarray*}
\end{proposition}
\textbf{Proof.}
Thanks to the monotonicity of $f$, the integrand in \eqref{velocity-representation} is nonnegative, that is
\begin{eqnarray*}
\frac{(f(y)-f(x))(x^2+3y^2+a^2)}{(x^2-y^2)((x-y)^2+a^2)((x+y)^2+a^2)}\geq0
\end{eqnarray*}
in either the regime $y<x$ or the regime $y>x$. Thus we can restrict the integral to the regime $0<y<x$, and obtain that
\begin{eqnarray*}
\begin{split}
H_af(x)
&\geq\frac{2a^2x}{\pi}\int\limits^x_0\frac{(f(y)-f(x))(x^2+3y^2+a^2)}{(x^2-y^2)((x-y)^2+a^2)((x+y)^2+a^2)}dy\\
&\geq\frac{2}{\pi}\cdot\frac{a^2}{x(x^2+a^2)}\int\limits^x_0(f(y)-f(x))dy,
\end{split}
\end{eqnarray*}
where in the last inequality we used the following formula: for any positive $y\neq x$,
\begin{eqnarray}\label{partial-derivative-y}
\begin{split}
&\frac{\partial}{\partial y}\Big(\frac{x^2+3y^2+a^2}{(x^2-y^2)((x-y)^2+a^2)((x+y)^2+a^2)}\Big)=\frac{6y(2x^2+2y^2+a^2)}{((x-y)^2+a^2)^2((x+y)^2+a^2)^2}\\
&\ \ \ \ \ \ \ \ \ \ \ \ \ \ \ \ \ \ \ \ \ \  \ \
+2a^2y\cdot\frac{2(3x^4+2x^2y^2+2a^2y^2+3y^4)+a^2(4x^2+a^2)}{(x^2-y^2)^2((x-y)^2+a^2)^2((x+y)^2+a^2)^2}>0.
\end{split}
\end{eqnarray}
This completes the proof of the proposition. \hfill $\square$
\begin{remark}
Another elementary proof but with a slightly inferior constant under the assumption of Proposition \ref{lower-bound-Ha} is as follows.

Indeed, by \eqref{definition-Ha}, integrating by parts along with the use of the identity
\begin{eqnarray}\label{derivative-identity-3.2}
\frac{d}{dx}\log\frac{|x|}{\sqrt{x^2+a^2}}=\frac{a^2}{x(x^2+a^2)},
\end{eqnarray}
and the even symmetry of $f$, it is not difficult to check that
\begin{eqnarray*}
\begin{split}
H_af(x)
&=\frac{1}{\pi}P.V.\int\limits_{\R}\log\Big(\frac{|x-y|}{\sqrt{(x-y)^2+a^2}}\Big)f'(y)dy\\
&=\frac{1}{2\pi}\int\limits_0^\infty\log\Big(\Big|\frac{x-y}{x+y}\Big|^2\frac{(x+y)^2+a^2}{(x-y)^2+a^2}\Big)f'(y)dy.
\end{split}
\end{eqnarray*}
Furthermore, by the elementary inequality $\log(1+t)\leq t$ for any $t>-1$ and $f'\leq0$, we obtain that
\begin{eqnarray*}
\begin{split}
H_af(x)
&\geq\frac{1}{2\pi}\int\limits_0^x\log\Big(1+\frac{-4a^2xy}{(x+y)^2((x-y)^2+a^2)}\Big)f'(y)dy\\
&\geq\frac{1}{\pi}\int\limits_0^x\frac{-2a^2xyf'(y)}{(x+y)^2((x-y)^2+a^2)}dy
\geq-\frac{1}{2\pi}\cdot\frac{a^2}{x(x^2+a^2)}\int\limits_0^xyf'(y)dy\\
&=\frac{1}{2\pi}\cdot\frac{a^2}{x(x^2+a^2)}\int\limits_0^x(f(y)-f(x))dy.
\end{split}
\end{eqnarray*}
\end{remark}
Proposition \ref{lower-bound-Ha} can now be used to establish our first weighted nonlinear integral inequality for the transform $H_a$.
\begin{proposition}\label{nonlinear-inequality-Ha-prop-3.2}
Let $f:\R\rightarrow\R$ be an even Schwartz function which is monotone decreasing on $(0,\infty)$. Then for any $a>0$ and $\sigma\in(-1,1)$, we have
\begin{eqnarray*}
-\int\limits^\infty_0\frac{H_a f(x) f'(x)}{x^{1+\sigma}}dx
\geq
C_{a,\sigma}\int\limits^\infty_0\frac{(f(x)-f(0))^2}{x^{2+\sigma}(x^2+a^2)}dx,
\end{eqnarray*}
where
\begin{eqnarray*}
C_{a,\sigma}=\frac{a^2}{\pi}(3+\sigma-2\sqrt{2+\sigma}).
\end{eqnarray*}
\end{proposition}
\textbf{Proof.} By Proposition \ref{lower-bound-Ha}, we have
\begin{eqnarray*}
\begin{split}
-\int\limits^\infty_0\frac{H_a f(x) f'(x)}{x^{1+\sigma}}dx
&\geq-\frac{2a^2}{\pi}\int\limits^\infty_0\frac{\int_0^x(f(y)-f(x))dy}{x^{2+\sigma}(x^2+a^2)}\cdot f'(x) dx\\
&=\frac{2a^2}{\pi}\int\limits^\infty_0\frac{\int_0^x(g(x)-g(y))dy}{x^{2+\sigma}(x^2+a^2)}\cdot g'(x) dx,
\end{split}
\end{eqnarray*}
where we denoted $g(x)=f(0)-f(x)$.

By Fubini's theorem and integration by parts, we can obtain that
\begin{eqnarray*}
\begin{split}
&\int\limits^\infty_0\frac{\int_0^x(g(x)-g(y))dy}{x^{2+\sigma}(x^2+a^2)}\cdot g'(x)dx\\
&=\frac{1}{2}\int\limits^\infty_0\Big(\int\limits_y^\infty \frac{\frac{d}{dx}((g(x)-g(y))^2)}{x^{2+\sigma}(x^2+a^2)}dx\Big)dy\\
&=\frac{2+\sigma}{2}\iint\limits_{x\geq y>0}\frac{(g(x)-g(y))^2}{x^{3+\sigma}(x^2+a^2)}dxdy
+\iint\limits_{x\geq y>0}\frac{(g(x)-g(y))^2}{x^{1+\sigma}(x^2+a^2)^2}dxdy\\
&\geq\frac{2+\sigma}{2}\iint\limits_{x\geq y>0}\frac{(g(x)-g(y))^2}{x^{3+\sigma}(x^2+a^2)}dxdy.
\end{split}
\end{eqnarray*}
Now using the elementary inequality
\begin{eqnarray}\label{young}
(a_1-a_2)^2\geq(1-\alpha)a^2_1+(1-\frac{1}{\alpha})a^2_2
\end{eqnarray}
 for any $a_1,a_2\in\R$ and any $0<\alpha<1$, we obtain
\begin{eqnarray*}
\begin{split}
\iint\limits_{x\geq y>0}\frac{(g(x)-g(y))^2}{x^{3+\sigma}(x^2+a^2)}dxdy
&\geq(1-\alpha)\iint\limits_{x\geq y>0}\frac{(g(x))^2dxdy}{x^{3+\sigma}(x^2+a^2)}+(1-\frac{1}{\alpha})\iint\limits_{x\geq y>0}\frac{(g(y))^2dxdy}{x^{3+\sigma}(x^2+a^2)}\\
&\geq(1-\alpha)\iint\limits_{x\geq y>0}\frac{(g(x))^2dxdy}{x^{3+\sigma}(x^2+a^2)}+(1-\frac{1}{\alpha})\iint\limits_{x\geq y>0}\frac{(g(y))^2dxdy}{x^{3+\sigma}(y^2+a^2)}\\
&=(1-\alpha)\int\limits_0^\infty\frac{(g(x))^2dx}{x^{2+\sigma}(x^2+a^2)}+(1-\frac{1}{\alpha})\frac{1}{2+\sigma}\int\limits_0^\infty\frac{(g(y))^2dy}{(y^2+a^2)y^{2+\sigma}}\\
&=\frac{3+\sigma}{2+\sigma}\int\limits_0^\infty\frac{(g(x))^2dx}{x^{2+\sigma}(x^2+a^2)}-(\alpha+\frac{1}{\alpha(2+\sigma)})\int\limits_0^\infty\frac{(g(x))^2dx}{x^{2+\sigma}(x^2+a^2)},
\end{split}
\end{eqnarray*}
which follows from the choice of $\alpha=(2+\sigma)^{-\frac{1}{2}}\in(0,1)$ that
\begin{eqnarray*}
\iint\limits_{x\geq y>0}\frac{(g(x)-g(y))^2}{x^{3+\sigma}(x^2+a^2)}dxdy
\geq(\frac{3+\sigma}{2+\sigma}-\frac{2}{\sqrt{2+\sigma}})\int\limits_0^\infty\frac{(g(x))^2dx}{x^{2+\sigma}(x^2+a^2)}.
\end{eqnarray*}
Therefore, we finally obtain that
\begin{eqnarray*}
\begin{split}
-\int\limits^\infty_0\frac{H_a f(x) f'(x)}{x^{1+\sigma}}dx
\geq\frac{a^2}{\pi}(3+\sigma-2\sqrt{2+\sigma})\int\limits_0^\infty\frac{(f(0)-f(x))^2dx}{x^{2+\sigma}(x^2+a^2)},
\end{split}
\end{eqnarray*}
which then finishes the proof of Proposition \ref{nonlinear-inequality-Ha-prop-3.2}.\hfill $\square$
\begin{remark}
Under the assumption in Proposition \ref{nonlinear-inequality-Ha-prop-3.2},
we observe that
for any finite $L>0$, it should hold that
\begin{eqnarray}\label{finite-interval}
-\int\limits^L_0\frac{H_a f(x) f'(x)}{x^{1+\sigma}}dx
\geq
C_{a,\sigma,L}\int\limits^L_0\frac{(f(x)-f(0))^2}{x^{2+\sigma}}dx.
\end{eqnarray}
Indeed, by Proposition \ref{lower-bound-Ha}, we can see that, for any $x\in(0,L]$
\begin{eqnarray*}
H_a f(x)\geq \frac{2a^2}{\pi(a^2+L^2)}\cdot\frac{1}{x}\int\limits_0^x(f(y)-f(x))dy.
\end{eqnarray*}
Then, by completely similar arguments as in the proof of Proposition \ref{nonlinear-inequality-Ha-prop-3.2}, but dropping some harmless boundary terms, we can get \eqref{finite-interval} with the constant
\begin{eqnarray*}
C_{a,\sigma,L}=\frac{a^2(3+\sigma-2\sqrt{2+\sigma})}{\pi(a^2+L^2)}.
\end{eqnarray*}
\end{remark}
Next we generalize the inequality obtained in Proposition \ref{nonlinear-inequality-Ha-prop-3.2} to more general Kiselev-type inequality for the transform $H_a$, whose prototype for the Hilbert transform $H$ was first established in \cite{[Kiselev]}.
\begin{proposition}\label{nonlinear-inequality-Kiselev}
Let $f:\R\rightarrow\R$ be an even Schwartz function which is monotone decreasing on $(0,\infty)$. Then for any $a>0$, $\sigma\in(-1,1)$ and $p\in(1,\infty)$, we have
\begin{eqnarray*}
-\int\limits^\infty_0\frac{H_a f(x) f'(x)(f(0)-f(x))^{p-1}}{x^{1+\sigma}}dx
\geq
C_{a,\sigma,p}\int\limits^\infty_0\frac{(f(0)-f(x))^{p+1}}{x^{2+\sigma}(x^2+a^2)}dx,
\end{eqnarray*}
where
\begin{eqnarray*}
C_{a,\sigma,p}=\frac{(1+\sigma)a^2}{(p+1)\pi}\Big(1-(\frac{2}{3+\sigma})^{\frac1p}\Big)^p.
\end{eqnarray*}
\end{proposition}
\textbf{Proof.} By Proposition \ref{lower-bound-Ha}, we have
\begin{eqnarray*}
\begin{split}
\rm LHS
&\geq-\frac{2a^2}{\pi}\int\limits^\infty_0\frac{\int_0^x(f(y)-f(x))dy}{x^{2+\sigma}(x^2+a^2)}\cdot f'(x)(f(0)-f(x))^{p-1} dx\\
&=\frac{2a^2}{\pi}\iint\limits_{x\geq y>0}\frac{g'(x)(g(x)-g(y))(g(x))^{p-1}}{x^{2+\sigma}(x^2+a^2)}dxdy,
\end{split}
\end{eqnarray*}
where we denoted $g(x)=f(0)-f(x)$.
Then $g$ is nonnegative and monotone increasing on $(0,\infty)$.
Then for any $x>y>0$, we have that $0\leq\frac{g(y)}{g(x)}\leq1$.
Note that for any $s\in[0,1]$, it holds that $1-s\geq (1-s)^{p}$. This in turn implies that for any $x>y>0$,
\begin{eqnarray*}
(g(x)-g(y))(g(x))^{p-1}\geq(g(x)-g(y))^p.
\end{eqnarray*}
It follows that
\begin{eqnarray*}
\begin{split}
&\iint\limits_{x\geq y>0}\frac{g'(x)(g(x)-g(y))(g(x))^{p-1}}{x^{2+\sigma}(x^2+a^2)}dxdy\\
&\geq\iint\limits_{x\geq y>0}\frac{g'(x)(g(x)-g(y))^p}{x^{2+\sigma}(x^2+a^2)}dxdy\\
&=\frac{1}{p+1}\int\limits^\infty_0\Big(\int_y^\infty \frac{\frac{d}{dx}((g(x)-g(y))^{p+1})}{x^{2+\sigma}(x^2+a^2)}dx\Big)dy\\
&=\frac{2+\sigma}{p+1}\iint\limits_{x\geq y>0}\frac{(g(x)-g(y))^{p+1}}{x^{3+\sigma}(x^2+a^2)}dxdy
+\frac{2}{p+1}\iint\limits_{x\geq y>0}\frac{(g(x)-g(y))^{p+1}}{x^{1+\sigma}(x^2+a^2)^2}dxdy\\
&\geq\frac{2+\sigma}{p+1}\iint\limits_{x\geq y>0}\frac{(g(x)-g(y))^{p+1}}{x^{3+\sigma}(x^2+a^2)}dxdy.
\end{split}
\end{eqnarray*}
Now note that for any $\beta>1$, we have that for $s\in[0,1]$,
\begin{eqnarray*}
(1-s)^{p+1}\geq c_\ast(1-\beta s^{p+1}),
\end{eqnarray*}
where
\begin{eqnarray*}
c_\ast=\Big(1-(\frac1\beta)^{\frac1p}\Big)^p.
\end{eqnarray*}
This in turn implies that for any $x>y>0$,
\begin{eqnarray*}
(g(x)-g(y))^{p+1}\geq c_\ast((g(x))^{p+1}-\beta(g(y))^{p+1}).
\end{eqnarray*}
Therefore, we further have that
\begin{eqnarray*}
\begin{split}
\iint\limits_{x\geq y>0}\frac{(g(x)-g(y))^{p+1}}{x^{3+\sigma}(x^2+a^2)}dxdy
&\geq c_\ast\iint\limits_{x\geq y>0}\frac{(g(x))^{p+1}dxdy}{x^{3+\sigma}(x^2+a^2)}-c_\ast\beta\iint\limits_{x\geq y>0}\frac{(g(y))^{p+1}dxdy}{x^{3+\sigma}(x^2+a^2)}\\
&\geq c_\ast\iint\limits_{x\geq y>0}\frac{(g(x))^{p+1}dxdy}{x^{3+\sigma}(x^2+a^2)}-c_\ast\beta\iint\limits_{x\geq y>0}\frac{(g(y))^{p+1}dxdy}{x^{3+\sigma}(y^2+a^2)}\\
&=c_\ast\int\limits_0^\infty\frac{(g(x))^{p+1}dx}{x^{2+\sigma}(x^2+a^2)}-\frac{c_\ast\beta}{2+\sigma}\int\limits_0^\infty\frac{(g(y))^2dy}{y^{2+\sigma}(y^2+a^2)}\\
&=c_\ast\Big(1-\frac{\beta}{2+\sigma}\Big)\int\limits_0^\infty\frac{(g(x))^{p+1}dx}{x^{2+\sigma}(x^2+a^2)},
\end{split}
\end{eqnarray*}
which follows from the choice of  $\beta=\frac{3+\sigma}{2}\in(1,\infty)$ that
\begin{eqnarray*}
\iint\limits_{x\geq y>0}\frac{(g(x)-g(y))^{p+1}}{x^{3+\sigma}(x^2+a^2)}dxdy
\geq \frac{(1+\sigma)c_\ast}{2(2+\sigma)}\int\limits_0^\infty\frac{(g(x))^{p+1}}{x^{2+\sigma}(x^2+a^2)}dx.
\end{eqnarray*}
Therefore, we finally obtain that
\begin{eqnarray*}
\begin{split}
{\rm LHS}
\geq\frac{(1+\sigma)a^2}{(p+1)\pi}\Big(1-(\frac{2}{3+\sigma})^{\frac1p}\Big)^p\int\limits_0^\infty\frac{(f(0)-f(x))^{p+1}}{x^{2+\sigma}(x^2+a^2)}dx,
\end{split}
\end{eqnarray*}
which then finishes the proof of Proposition \ref{nonlinear-inequality-Kiselev}.   \hfill $\square$

Proposition \ref{nonlinear-inequality-Ha-prop-3.2} can be used to establish the finite time blow-up. Consider
\begin{equation}\label{KS-}
\left\{\ba
&\partial_t\rho-gH_a\rho\partial_x\rho= 0, ~(x,t)\in \R\times\R_+,\\
&\rho(x,0)=\rho_{0}(x),~x\in\R. \ea\ \right.
\end{equation}
It should be remarked that here we consider the equation with a minus sign for the nonlinear term rather than the a prior more physically meaningful with a  plus to simplify our presentation. It is clear that due to the properties of the transform $H_a$ the transformation $\rho\rightarrow-\rho$ transforms one equation into the other.
It is clear that Proposition \ref{local-well-posedness} and Lemmas \ref{nondecreasing-property}-\ref{maximum-principle} also hold true for the model \eqref{KS-}.

Our finite time blow-up result can be stated as
\begin{theorem}\label{singularity-formation-Li-Rodrigo-monotone}
Fix $a,g>0$ and $\sigma\in(0,1)$.  Let the initial data $\rho_0$ be an even Schwartz function on $\R$ and monotone decreasing on $(0,\infty)$.
If
\begin{eqnarray*}
\int\limits^\infty_{0}\frac{\rho_0(0)-\rho_0(x)}{x^{1+\sigma}}dx>\frac{2\sqrt{2}}{\sigma}\|\rho_0\|_{L^\infty},
\end{eqnarray*}
Then the smooth solution $\rho$ to \eqref{KS-} blows up in finite time.
\end{theorem}
\textbf{Proof.} Suppose the initial data $\rho_0\in\mathcal{S}(\R)$ is even and monotone decreasing on $(0,\infty)$,
 and let $\rho(x,t)$ denote the corresponding unique local solution to \eqref{KS-}. Now, assume for contradiction that $\rho$ exists for all time. Then $\rho$ is even and monotone decreasing on $(0,\infty)$ for all time $t>0$.
We define the quantity $J(t)$ as follows,
\begin{eqnarray*}
J(t)=\int\limits_0^\infty\frac{\rho(0,t)-\rho(x,t)}{x^{1+\sigma}}dx,
\end{eqnarray*}
where $\sigma\in(0,1)$ is arbitrarily fixed. We compute and invoke Proposition \ref{nonlinear-inequality-Ha-prop-3.2} to derive that
\begin{eqnarray}\label{Jt-ODE-monotone}
\begin{split}
\frac{d}{dt}J(t)
&=-g\int\limits^\infty_0\frac{H_a \rho(x,t)\partial_x \rho(x,t)}{x^{1+\sigma}}dx\\
&\geq gC_{a,\sigma}\int\limits^\infty_0\frac{(\rho(0,t)-\rho(x,t))^2}{x^{2+\sigma}(x^2+a^2)}dx.
\end{split}
\end{eqnarray}
By the Cauchy-Schwarz inequality and Lemma \ref{maximum-principle}, $J(t)$ can be bounded by
\begin{eqnarray*}
\begin{split}
0\leq J(t)
&\leq\Big(\int\limits^1_0\frac{(\rho(0,t)-\rho(x,t))^2}{x^{2+\sigma}(x^2+a^2)}dx\Big)^{\frac12}\Big(\int\limits^1_0\frac{x^2+a^2}{x^\sigma}dx\Big)^{\frac12}+2\|\rho_0\|_{L^\infty}\int\limits_1^\infty\frac{dx}{x^{1+\sigma}}\\
&=\Big(\frac{1}{3-\sigma}+\frac{a^2}{1-\sigma}\Big)^{\frac12}\Big(\int\limits^1_0\frac{(\rho(0,t)-\rho(x,t))^2}{x^{2+\sigma}(x^2+a^2)}dx\Big)^{\frac12}+\frac{2}{\sigma}\|\rho_0\|_{L^\infty},
\end{split}
\end{eqnarray*}
which along with \eqref{Jt-ODE-monotone} implies that
\begin{eqnarray*}
\begin{split}
\frac{d}{dt}J(t)
\geq\frac{g(1-\sigma)(3-\sigma)C_{a,\sigma}}{2(1-\sigma+(3-\sigma)a^2)}(J(t))^2
-\frac{4g(1-\sigma)(3-\sigma)C_{a,\sigma}}{(1-\sigma+(3-\sigma)a^2)\sigma^2}\|\rho_0\|^2_{L^\infty}.
\end{split}
\end{eqnarray*}
From this differential inequality, we can deduce that $J(t)$ will blow up in finite time, leading to a contradiction, provided
\begin{eqnarray*}
J(0)=\int\limits_0^\infty\frac{\rho_0(0)-\rho_0(x)}{x^{1+\sigma}}dx>\frac{2\sqrt{2}}{\sigma}\|\rho_0\|_{L^\infty}.
\end{eqnarray*}
We then have proved the theorem.    \hfill $\square$
\subsection{Non-monotone decaying case}
In this subsection, we give another proof of the finite blow-up which works for only even symmetric (not necessarily monotone decaying) initial data. We begin with another proof of the inequality in Proposition \ref{nonlinear-inequality-Ha-prop-3.2} with $\sigma=0$.
\begin{proposition}\label{1-nonlinear-inequality-nodecay}
Let $f:\R\rightarrow\R$ be an even Schwartz function. Then for every $a>0$, we have
\begin{eqnarray*}
-\int\limits^\infty_0\frac{H_a f(x) f'(x)}{x}dx
\geq
\frac{a^2}{\pi}
\int\limits^\infty_0\frac{(f(0)-f(x))^2}{x^2}\cdot\frac{3x^2+a^2}{(x^2+a^2)^2}dx.
\end{eqnarray*}
\end{proposition}
\textbf{Proof.}
By \eqref{velocity-representation} and integration by parts, we derive that
\begin{align*}
&-\int\limits^\infty_0\frac{H_a f(x) f'(x)}{x}dx\\
&=\frac{2a^2}{\pi}\int\limits^\infty_0f'(x)\Big(\int\limits^\infty_0\frac{(f(x)-f(y))(x^2+3y^2+a^2)}{(x^2-y^2)((x-y)^2+a^2)((x+y)^2+a^2)}dy\Big)dx\\
&=\frac{a^2}{\pi}\lim_{\epsilon\rightarrow0}\int\limits^\infty_0\Big(\int\limits^\infty_0\frac{(x^2+3y^2+a^2)\frac{\partial}{\partial x}((f(x)-f(y))^2)}{(x^2-y^2+i\epsilon)((x-y)^2+a^2)((x+y)^2+a^2)}dx\Big)dy\\
&=\frac{a^2}{\pi}\int\limits_0^\infty\frac{(f(0)-f(y))^2}{y^2}\cdot\frac{3y^2+a^2}{(y^2+a^2)^2}dy\\
&\ \ \
-\frac{a^2}{\pi}\int\limits^\infty_0\int\limits^\infty_0(f(x)-f(y))^2\frac{\partial}{\partial  x}\Big(\frac{x^2+3y^2+a^2}{(x^2-y^2)((x-y)^2+a^2)((x+y)^2+a^2)}\Big)dxdy\\
&\geq\frac{a^2}{\pi}\int_0^\infty\frac{(f(0)-f(y))^2}{y^2}\cdot\frac{3y^2+a^2}{(y^2+a^2)^2}dy,
\end{align*}
where in the last inequality we used the following formula: for any positive $x\neq y$
\begin{eqnarray}\label{partial-derivative-x}
\begin{split}
&\frac{\partial}{\partial x}\Big(\frac{x^2+3y^2+a^2}{(x^2-y^2)((x-y)^2+a^2)((x+y)^2+a^2)}\Big)=-\frac{2x(2x^2+10y^2+3a^2)}{((x-y)^2+a^2)^2((x+y)^2+a^2)^2}\\
&\ \ \ \ \ \ \ \ \ \ \ \ \ \ \ \ \ \ \ \ \ \ \ \ \ \ \ \ \
-2a^2x\cdot\frac{2(x^4+6x^2y^2+2a^2x^2+y^4)+a^2(4y^2+a^2)}{(x^2-y^2)^2((x-y)^2+a^2)^2((x+y)^2+a^2)^2}<0.
\end{split}
\end{eqnarray}
Then we have proved the proposition.  \hfill $\square$
\par
Proposition \ref{1-nonlinear-inequality-nodecay} is not directly useful for establishing finite time blow-ups since it involves a non-integrable weight $\frac{1}{x}$. The next proposition fixes this issue.
\begin{proposition}\label{nonlinear-weighted-inequality-exponential-prop-3.6}
Let $f:\R\rightarrow\R$ be an even Schwartz function. Then, for every $a>0$, we have
\begin{eqnarray*}
-\int\limits^\infty_0\frac{H_af(x)f'(x)}{x}e^{-x}dx
\geq \frac{a^2}{2\pi}
\int\limits\limits^\infty_0\frac{(f(0)-f(x))^2(3x^2+a^2)}{x^2(x^2+a^2)^2}dx
-\frac{2}{\pi}(4719+\frac{3}{a^2})\|f\|^2_{L^\infty}.
\end{eqnarray*}
\end{proposition}
\textbf{Proof.} By using the same integration by parts argument as in Proposition \ref{1-nonlinear-inequality-nodecay} and the formula \eqref{partial-derivative-x}, we can get that
\begin{eqnarray*}
\begin{split}
&-\int\limits^\infty_0\frac{H_a f(x) f'(x)}{x}e^{-x}dx
=\frac{a^2}{\pi}\int\limits_0^\infty\frac{(f(0)-f(y))^2}{y^2}\cdot\frac{3y^2+a^2}{(y^2+a^2)^2}dy\\
&\ \ \
-\frac{a^2}{\pi}\int\limits^\infty_0\int\limits^\infty_0(f(x)-f(y))^2e^{-x}\frac{\partial}{\partial  x}\Big(\frac{x^2+3y^2+a^2}{(x^2-y^2)((x-y)^2+a^2)((x+y)^2+a^2)}\Big)dxdy\\
&\ \ \ \ \ \ \ \ \ \ \
+\frac{a^2}{\pi}\int\limits_0^\infty\int\limits_0^\infty \frac{(f(x)-f(y))^2(x^2+3y^2+a^2)e^{-x}}{(x^2-y^2)((x-y)^2+a^2)((x+y)^2+a^2)}dxdy\\
&\geq\frac{a^2}{\pi}\int\limits_0^\infty\frac{(f(0)-f(y))^2}{y^2}\cdot\frac{3y^2+a^2}{(y^2+a^2)^2}dy\\
&\ \ \
+\frac{a^2}{2\pi}\int\limits_0^\infty\int\limits_0^\infty \frac{(f(x)-f(y))^2((x^2+3y^2+a^2)e^{-x}-(y^2+3x^2+a^2)e^{-y})}{(x^2-y^2)((x-y)^2+a^2)((x+y)^2+a^2)}dxdy.
\end{split}
\end{eqnarray*}
Note that
\begin{eqnarray*}
\begin{split}
&\Big|\iint\limits_{0<\frac{x}{2}<y<2x} \frac{(f(x)-f(y))^2((x^2+3y^2+a^2)e^{-x}-(y^2+3x^2+a^2)e^{-y})}{(x^2-y^2)((x-y)^2+a^2)((x+y)^2+a^2)}dxdy\Big|\\
&\leq\Big|\iint\limits_{0<\frac{x}{2}<y<2x} \frac{(f(x)-f(y))^2(x^2+3y^2+a^2)(e^{-x}-e^{-y})}{(x^2-y^2)((x-y)^2+a^2)((x+y)^2+a^2)}dxdy\Big|\\
&\ \ \ \ \
+2\iint\limits_{0<\frac{x}{2}<y<2x} \frac{(f(x)-f(y))^2e^{-y}}{((x-y)^2+a^2)((x+y)^2+a^2)}dxdy\\
&\leq\frac{12}{a^2}\|f\|^2_{L^\infty}\iint\limits_{0<\frac{x}{2}<y<2x} \frac{|e^{-x}-e^{-y}|}{|x^2-y^2|}dxdy
+\frac{8}{a^4}\|f\|^2_{L^\infty}\iint\limits_{0<\frac{x}{2}<y<2x}e^{-y}dxdy\\
&\leq\frac{12}{a^2}\|f\|^2_{L^\infty}
+\frac{12}{a^4}\|f\|^2_{L^\infty}=\frac{12}{a^2}(1+\frac{1}{a^2})\|f\|^2_{L^\infty}.
\end{split}
\end{eqnarray*}
We denote the weight
\begin{eqnarray*}
\omega_a(x):=\frac{3x^2+a^2}{(x^2+a^2)^2}.
\end{eqnarray*}
Then it is not difficult to check that, for every $x>0$,
\begin{eqnarray}\label{weight-upper-bound}
0<\omega_a(x)\leq\frac{9}{8a^2}.
\end{eqnarray}
Note that
\begin{equation*}
\left\{\ba
&(y^2+3x^2+a^2)e^{-y}>(x^2+3y^2+a^2)e^{-x}, ~{\rm if}~ 0<y<x,\\
&(x^2+3y^2+a^2)e^{-x}>(y^2+3x^2+a^2)e^{-y}, ~{\rm if}~ 0<x<y.\ea\ \right.
\end{equation*}
It follows  that
\begin{eqnarray*}
\begin{split}
&\Big|\iint\limits_{(y\leq\frac{x}{2})\cup (y\geq2x)} \frac{(f(x)-f(y))^2((x^2+3y^2+a^2)e^{-x}-(y^2+3x^2+a^2)e^{-y})}{(x^2-y^2)((x-y)^2+a^2)((x+y)^2+a^2)}dxdy\Big|\\
&\leq\iint\limits_{0<y\leq\frac{x}{2}}\frac{(f(x)-f(y))^2(y^2+3x^2+a^2)}{(x^2-y^2)((x-y)^2+a^2)((x+y)^2+a^2)}\cdot e^{-y}dxdy\\
&\ \ \ \ \ \
+\iint\limits_{y\geq2x>0}\frac{(f(x)-f(y))^2(x^2+3y^2+a^2)}{(y^2-x^2)((x-y)^2+a^2)((x+y)^2+a^2)}\cdot e^{-x}dxdy\\
&=2\iint\limits_{0<y\leq\frac{x}{2}}\frac{(f(x)-f(y))^2(y^2+3x^2+a^2)}{(x^2-y^2)((x-y)^2+a^2)((x+y)^2+a^2)}\cdot e^{-y}dxdy\\
&\leq32\iint\limits_{0<y\leq\frac{x}{2}}\frac{(f(x)-f(y))^2e^{-y}}{x^2(x^2+4a^2)}dxdy
\leq32\iint\limits_{0<y\leq\frac{x}{2}} \frac{(f(x)-f(y))^2}{x^2}\omega_a(x)e^{-y}dxdy\\
&\leq64\underbrace{\int\limits^\infty_0 \frac{(f(x)-f(0))^2}{x^2}\omega_a(x)(1-e^{-\frac{x}{2}})dx}_I+64\underbrace{\int\limits^\infty_0\frac{(f(0)-f(y))^2}{e^y}\int\limits^\infty_{2y}\frac{\omega_a(x)}{x^2}dxdy}_{II},
\end{split}
\end{eqnarray*}
where in the last inequality we used the simple fact $(\xi+\eta)^2\leq2(\xi^2+\eta^2)$ for any $\xi,\eta\in\mathbb{R}$.

Furthermore, by \eqref{weight-upper-bound} and Cauchy-Schwarz inequality, we have
\begin{eqnarray*}
\begin{split}
I
&\leq\|f\|_{L^\infty}\int\limits^1_0 \frac{|f(x)-f(0)|}{x}\omega_a(x)dx
+\frac{9}{2a^2}\|f\|^2_{L^\infty}\int\limits^\infty_1 \frac{dx}{x^2}\\
&\leq\|f\|_{L^\infty}\Big(\int\limits^1_0 \frac{(f(x)-f(0))^2}{x^2}\omega_a(x)dx\Big)^{\frac12}\Big(\int\limits^1_0\omega_a(x)dx\Big)^{\frac12}
+\frac{9}{2a^2}\|f\|^2_{L^\infty}\\
&\leq\frac{3}{2\sqrt{2}a}\|f\|_{L^\infty}\Big(\int\limits^1_0 \frac{(f(x)-f(0))^2}{x^2}\omega_a(x)dx\Big)^{\frac12}
+\frac{9}{2a^2}\|f\|^2_{L^\infty}
\end{split}
\end{eqnarray*}
and
\begin{eqnarray*}
\begin{split}
II
&\leq\frac{3}{2}\int\limits^1_0\frac{(f(0)-f(y))^2}{y(y^2+a^2)}dy+\frac{9}{2a^2}\|f\|^2_{L^\infty}\int\limits^\infty_1\frac{1}{y}\Big(\int\limits^\infty_{2y}\frac{1}{x^2}dx\Big)dy\\
&\leq3\|f\|_{L^\infty}\int\limits^1_0\frac{|f(0)-f(y)|}{y}\omega_a(y)dy+\frac{9}{4a^2}\|f\|^2_{L^\infty}\\
&\leq\frac{9}{2\sqrt{2}a}\|f\|_{L^\infty}\Big(\int\limits^1_0 \frac{(f(x)-f(0))^2}{x^2}\omega_a(x)dx\Big)^{\frac12}+\frac{9}{4a^2}\|f\|^2_{L^\infty}.
\end{split}
\end{eqnarray*}
It follows  that
\begin{eqnarray*}
\begin{split}
&\frac{a^2}{2\pi}\Big|\iint\limits_{(y\leq\frac{x}{2})\cup (y\geq2x)}\frac{(f(x)-f(y))^2((x^2+3y^2+a^2)e^{-x}-(y^2+3x^2+a^2)e^{-y})}{(x^2-y^2)((x-y)^2+a^2)((x+y)^2+a^2)} dxdy\Big|\\
&\leq\frac{96\sqrt{2}a}{\pi}\|f\|_{L^\infty}\Big(\int\limits^1_0 \frac{(f(x)-f(0))^2}{x^2}\omega_a(x)dx\Big)^{\frac12}+\frac{216}{\pi}\|f\|^2_{L^\infty}.
\end{split}
\end{eqnarray*}
Finally, by Young's inequality, we have that
\begin{eqnarray*}
\begin{split}
&\frac{a^2}{2\pi}\Big|\int\limits_0^\infty\int\limits_0^\infty \frac{(f(x)-f(y))^2((x^2+3y^2+a^2)e^{-x}-(y^2+3x^2+a^2)e^{-y})}{(x^2-y^2)((x-y)^2+a^2)((x+y)^2+a^2)}dxdy\Big|\\
&\leq\frac{96\sqrt{2}a}{\pi}\|f\|_{L^\infty}\Big(\int\limits^1_0 \frac{(f(x)-f(0))^2}{x^2}\omega_a(x)dx\Big)^{\frac12}+\frac{2}{\pi}(111+\frac{3}{a^2})\|f\|^2_{L^\infty}\\
&\leq\frac{a^2}{2\pi}\int\limits^1_0 \frac{(f(x)-f(0))^2}{x^2}\omega_a(x)dx+\frac{2}{\pi}(4719+\frac{3}{a^2})\|f\|^2_{L^\infty}.
\end{split}
\end{eqnarray*}
We then conclude the proof of this proposition.   \hfill $\square$

Proposition \ref{nonlinear-weighted-inequality-exponential-prop-3.6} can be used to establish the blow-up of the solutions to \eqref{KS} for the initial data which is not necessarily monotone on $(0,\infty)$. The main result of this aspect is stated as
\begin{theorem}\label{singularity-formation-Li-Rodrigo-non-monotone}
Let the initial data $\rho_0$ be an even Schwartz function. There exists a constant $A_a>0$ depending only on $a$ such that if
\begin{eqnarray*}
\int\limits^\infty_{0}\frac{\rho_0(x)-\rho_0(0)}{x}e^{-x}dx\geq A_a\|\rho_0\|_{L^\infty},
\end{eqnarray*}
Then the smooth solution $\rho$ to \eqref{KS} blows up in finite time.
\end{theorem}
\textbf{Proof.} Define the quantity $\widetilde{J}(t)$ as follows,
\begin{eqnarray*}
\widetilde{J}(t)=\int\limits_0^\infty\frac{\rho(x,t)-\rho(0,t)}{x}e^{-x}dx.
\end{eqnarray*}
By using Proposition \ref{nonlinear-weighted-inequality-exponential-prop-3.6} and Lemma \ref{maximum-principle}, we compute
\begin{eqnarray}\label{Jt-ODE}
\begin{split}
\frac{d}{dt}\widetilde{J}(t)
&=-g\int\limits^\infty_0\frac{H_a \rho(x,t)\partial_x \rho(x,t)}{x}e^{-x}dx\\
&\geq \frac{a^2g}{2\pi}\int\limits^\infty_0\frac{(\rho(0,t)-\rho(x,t))^2(3x^2+a^2)}{x^2(x^2+a^2)^2}dx
-\frac{2g}{\pi}(4719+\frac{3}{a^2}) \|\rho_0\|^2_{L^\infty}.
\end{split}
\end{eqnarray}
By the Cauchy-Schwarz inequality, $\widetilde{J}(t)$ can be bounded by
\begin{eqnarray*}
\begin{split}
|\widetilde{J}(t)|
&\leq\Big(\int\limits^\infty_0\frac{(\rho(0,t)-\rho(x,t))^2(3x^2+a^2)}{x^{2}(x^2+a^2)^2}dx\Big)^{\frac12}
\Big(\int\limits^\infty_0\frac{(x^2+a^2)^2}{3x^2+a^2}e^{-2x}dx\Big)^{\frac12}\\
&\leq\Big(\frac{1}{4}+\frac{a^2}{2}\Big)^{\frac12}
\Big(\int\limits^\infty_0\frac{(\rho(0,t)-\rho(x,t))^2(3x^2+a^2)}{x^2(x^2+a^2)}dx\Big)^{\frac12},
\end{split}
\end{eqnarray*}
which along with \eqref{Jt-ODE} implies that
\begin{eqnarray*}
\begin{split}
\frac{d}{dt}\widetilde{J}(t)
\geq c_1(\widetilde{J}(t))^2-c_2,~{\rm with}~c_1:=\frac{2a^2g}{(1+2a^2)\pi},~c_2:=\frac{2g}{\pi}(4719+\frac{3}{a^2}) \|\rho_0\|^2_{L^\infty}.
\end{split}
\end{eqnarray*}
With this differential inequality in hand, following the arguments in \cite{[Jarrin-Vegara-Hermosilla]}, we conclude that the functional $\widetilde{J}(t)$ blows-up at the finite time
\begin{eqnarray*}
t_\ast=\frac{1}{2\sqrt{c_1c_2}}\log\frac{\sqrt{c_1}\widetilde{J}(0)+\sqrt{c_2}}{\sqrt{c_1}\widetilde{J}(0)-\sqrt{c_2}}>0,
\end{eqnarray*}
provided
\begin{eqnarray*}
\widetilde{J}(0)=
\int\limits_0^\infty\frac{\rho_0(x)-\rho_0(0)}{x}e^{-x}dx
>\sqrt{\frac{c_2}{c_1}}=\sqrt{(3+4719a^2)(1+2a^2)}\|\rho_0\|_{L^\infty}.
\end{eqnarray*}
We then have proved the theorem.   \hfill $\square$
\begin{remark}
We note that our results can be extended to equations with fractional dissipations:
\begin{eqnarray}\label{dissipative-IPM}
\partial_t\rho+gH_a\rho\partial_x\rho+(-\Delta)^{\frac{\gamma}{2}}\rho=0.
\end{eqnarray}
By using Proposition \ref{nonlinear-weighted-inequality-exponential-prop-3.6} and following the proof of Theorem 4.5 in \cite{[Li-Rodrigo20]}, one can show the following finite-time blow-up for \eqref{dissipative-IPM} with $\gamma\in(0,\frac12)$. Let $\gamma\in(0,\frac12)$. Let the initial data $\rho_0$ be an even Schwartz function.
There exists a constant $A_a>0$ depending only on $a$ such that if
\begin{eqnarray*}
\int\limits^\infty_{0}\frac{\rho_0(x)-\rho_0(0)}{x}e^{-x}dx\geq A_a\|\rho_0\|_{L^\infty},
\end{eqnarray*}
Then the smooth solution $\rho$ to \eqref{dissipative-IPM} blows up in finite time.
\end{remark}
\section{Singularity formation via telescoping sums}
In this section we present another proof of finite time singularity formation for smooth solutions to \eqref{KS-} (Theorem \ref{singularity-formation-telescoping-sums}). The proof is based on a local in space upper bound for the nonlinearity, which is established in Lemma \ref{upper-bound-nonlinear-vicol}.
\begin{lemma}\label{upper-bound-Ha}
Let $f:\R\rightarrow\R$ be even and monotone decreasing on $(0,\infty)$. Then for any fixed $a>0$ and $x_1>x_2>0$, we have
\begin{eqnarray*}
H_a f(x_2)\geq \frac{1}{\pi}(f(x_1)-f(x_2))\log\Big(\frac{x_1-x_2}{x_1+x_2}\sqrt{\frac{(x_1+x_2)^2+a^2}{(x_1-x_2)^2+a^2}}\Big).
\end{eqnarray*}
\end{lemma}
\textbf{Proof.} By \eqref{velocity-representation} and the monotonicity of $f$, we derive the lower bound as follows:
\begin{eqnarray*}
\begin{split}
H_af(x_2)
&\geq\frac{2a^2x_2}{\pi}\int\limits^\infty_{x_1}\frac{(f(y)-f(x_2))(x^2_2+3y^2+a^2)}{(x^2_2-y^2)((x_2-y)^2+a^2)((x_2+y)^2+a^2)}dy\\
&\geq\frac{1}{\pi}(f(x_1)-f(x_2))\int\limits^\infty_{x_1}\frac{2a^2x_2(x^2_2+3y^2+a^2)}{(x^2_2-y^2)((x_2-y)^2+a^2)((x_2+y)^2+a^2)}dy\\
&=\frac{1}{\pi}(f(x_1)-f(x_2))\log\Big(\frac{x_1-x_2}{x_1+x_2}\sqrt{\frac{(x_1+x_2)^2+a^2}{(x_1-x_2)^2+a^2}}\Big),
\end{split}
\end{eqnarray*}
which is the desired lower bound.   \hfill $\square$
\par
The next lemma yields a upper bound for the nonlinear term in \eqref{KS}, which is local in nature.
\begin{lemma}\label{upper-bound-nonlinear-vicol}
Let $f:\R\rightarrow\R$ be even and monotone decreasing on $(0,\infty)$. Then for any fixed $a>0$ and $x_1>x_2>0$, we have
\begin{eqnarray*}
\int\limits_{x_2}^{x_1}f'(x)H_a f(x)dx\leq \frac{1}{4\pi}(f(x_1)-f(x_2))^2\log\Big(\frac{x_1-x_2}{x_1+x_2}\sqrt{\frac{(x_1+x_2)^2+a^2}{(x_1-x_2)^2+a^2}}\Big)
\end{eqnarray*}
\end{lemma}
\textbf{Proof.}
Choose a point $x_3\in(x_2,x_1)$ such that $f(x_3)=\frac{1}{2}(f(x_1)+f(x_2))$. Without loss of generality we can assume
$x_3\leq\frac{1}{2}(x_1+x_2)$, as the other case can be treated similarly. By Lemma \ref{upper-bound-Ha} and the monotonicity of $f$, we obtain that, for any $x\in(x_2,x_3)$
\begin{eqnarray*}
\begin{split}
H_a f(x)
&\geq\frac{1}{\pi}(f(x_1)-f(x))\log\Big(\frac{x_1-x}{x_1+x}\sqrt{\frac{(x_1+x)^2+a^2}{(x_1-x)^2+a^2}}\Big)\\
&\geq\frac{1}{\pi}(f(x_1)-f(x_3))\log\Big(\frac{x_1-x}{x_1+x}\sqrt{\frac{(x_1+x)^2+a^2}{(x_1-x)^2+a^2}}\Big)\\
&\geq\frac{1}{2\pi}(f(x_1)-f(x_2))\log\Big(\frac{x_1-x_2}{x_1+x_2}\sqrt{\frac{(x_1+x_2)^2+a^2}{(x_1-x_2)^2+a^2}}\Big),
\end{split}
\end{eqnarray*}
where we used the simple fact that, for any $x\in(x_2,x_1)$
\begin{eqnarray*}
\frac{d}{dx}\log\Big(\frac{x_1-x}{x_1+x}\sqrt{\frac{(x_1+x)^2+a^2}{(x_1-x)^2+a^2}}\Big)
=\frac{2a^2x_1(x^2_1+3x^2+a^2)}{(x^2-x^2_1)((x-x_1)^2+a^2)((x+x_1)^2+a^2)}<0.
\end{eqnarray*}
Note that $H_a f(x)f'(x)\leq0$ by our assumptions on $f$. It follows from the choice of $x_3$ that
\begin{eqnarray*}
\begin{split}
\int\limits_{x_2}^{x_1}H_a f(x)f'(x)dx
&\leq\int\limits_{x_2}^{x_3}H_a f(x)f'(x)dx\\
&\leq\frac{1}{2\pi}(f(x_1)-f(x_2))\log\Big(\frac{x_1-x_2}{x_1+x_2}\sqrt{\frac{(x_1+x_2)^2+a^2}{(x_1-x_2)^2+a^2}}\Big)\int\limits_{x_2}^{x_3}f'(x)dx\\
&=\frac{1}{4\pi}(f(x_1)-f(x_2))^2\log\Big(\frac{x_1-x_2}{x_1+x_2}\sqrt{\frac{(x_1+x_2)^2+a^2}{(x_1-x_2)^2+a^2}}\Big),
\end{split}
\end{eqnarray*}
which is the wished upper bound. \hfill$\square$
\begin{remark}
Note that for any Schwartz function $f$, by Fubini's theorem and  integration by parts, we can derive that
\begin{eqnarray*}
\begin{split}
\int\limits_{\R}f'(x)H_a f(x)dx
&=\frac{a^2}{\pi}\int\limits_{\R}f'(x)\int\limits_{\R}\frac{f(y)-f(x)}{(x-y)((x-y)^2+a^2)}dydx\\
&=-\frac{a^2}{2\pi}\lim_{\epsilon\rightarrow0}\int\limits_{\R}\int\limits_{\R}\frac{\frac{\partial}{\partial x}(f(x)-f(y))^2}{(x-y+i\epsilon)((x-y)^2+a^2)}dxdy\\
&=-\frac{a^2}{2\pi}\iint\limits_{\R\times\R}\frac{(f(x)-f(y))^2(3(x-y)^2+a^2)}{(x-y)^2((x-y)^2+a^2)^2}dxdy.
\end{split}
\end{eqnarray*}
Then Lemma \ref{upper-bound-nonlinear-vicol} can be viewed as a local version of this above identity. Note that in the proof of Lemma \ref{upper-bound-nonlinear-vicol} we did not integrate by parts.
Also, from the above identity, we can see that the integral of the solution $\rho$ to \eqref{KS} must be increasing.
\end{remark}
With Lemma \ref{upper-bound-nonlinear-vicol} in hand, we can prove the finite time singularity by the method of telescoping sums. Motivated by \cite{[Silvestre-Vicol]}, the idea of the proof is to use a weighted version of Lemma \ref{upper-bound-nonlinear-vicol} in a dyadic fashion.

\begin{theorem}\label{singularity-formation-telescoping-sums}
Let $\rho_0$ be an even Schwartz function which is monotone decreasing on $(0,\infty)$. Then the Cauchy problem for \eqref{KS-} does not have a global in time smooth solution.
\end{theorem}
\textbf{Proof.}
Consider the continuous function $\eta:(0,\infty)\rightarrow(0,\infty)$ defined by
\begin{eqnarray*}
\eta(x)
=
\begin{cases}
x^{-\alpha},
&  \mbox{if $0<x<1,$ }\\
x^{-4-\alpha},
& \mbox{if $x\geq1,$ } \\
\end{cases}
\end{eqnarray*}
for some fixed $\alpha\in(0,1)$. Then $\eta\in L ^1(\mathbb{R}_+)$ is non-increasing on $(0,\infty)$.
Define the weighted integral
\begin{eqnarray*}
F(t)\equiv\int\limits_0^\infty \eta(x)(\rho(0,t)-\rho(x,t))dx\geq0.
\end{eqnarray*}
Then, in view of Lemma \ref{maximum-principle}, we have that
\begin{eqnarray}\label{bounded-F}
|F(t)|\leq2\|\rho_0\|_{L^\infty}\|\eta\|_{L^1(\mathbb{R}_+)}=\frac{8\|\rho_0\|_{L^\infty}}{(1-\alpha)(3+\alpha)}.
\end{eqnarray}
On the other hand, by \eqref{KS-} and the fact $H_a\rho(0,t)=0$, we compute
\begin{eqnarray}\label{ODE-F}
\begin{split}
\frac{d}{dt}F(t)
&=\int\limits_0^\infty\eta(x)\partial_t(\rho(0,t)-\rho(x,t))dx\\
&=-g\int\limits_0^\infty\eta(x)H_a\rho(x,t)\partial_x\rho(x,t)dx+gH_a\rho(0,t)\partial_x\rho(0,t)\int\limits_0^\infty\eta(x)dx\\
&=-g\sum_{k\in\mathbb{Z}}\int\limits_{2^k}^{2^{k+1}}\eta(x)H_a\rho(x,t)\partial_x\rho(x,t)dx.
\end{split}
\end{eqnarray}
By Lemma \ref{upper-bound-nonlinear-vicol} and the monotonicity of $\eta$ and $\rho$, we further derive that
\begin{eqnarray*}
\begin{split}
&-\int\limits_{2^k}^{2^{k+1}}\eta(x)H_a\rho(x,t)\partial_x\rho(x,t)dx\\
&\geq-\eta(2^{k+1})\int\limits_{2^k}^{2^{k+1}}H_a\rho(x,t)\partial_x\rho(x,t)dx\\
&\geq\frac{\eta(2^{k+1})}{4\pi}\log\Big(3\sqrt{\frac{2^{2k}+a^2}{9\cdot2^{2k}+a^2}}\Big)(\rho(2^{k+1},t)-\rho(2^k,t))^2,
\end{split}
\end{eqnarray*}
which along with \eqref{ODE-F} implies that
\begin{eqnarray}\label{ODE-F-inequality}
\frac{d}{dt}F(t)\geq\frac{g}{4\pi}\sum_{k\in\mathbb{Z}}\eta(2^{k+1})\log\Big(3\sqrt{\frac{2^{2k}+a^2}{9\cdot2^{2k}+a^2}}\Big)(\rho(2^{k+1},t)-\rho(2^k,t))^2.
\end{eqnarray}
In order to obtain a lower bound of the infinite series (in \eqref{ODE-F-inequality}) in terms of $F(t)$, we introduce the function $\varphi$ defined by
\begin{eqnarray*}
\varphi(x)
=
\begin{cases}
\frac{1}{3+\alpha}+\frac{1}{1-\alpha}(1-x^{1-\alpha}),
&  \mbox{if $0<x<1,$ }\\
\frac{1}{3+\alpha}x^{-3-\alpha},
& \mbox{if $x\geq1.$ } \\
\end{cases}
\end{eqnarray*}
It is clear that
\begin{eqnarray*}
\eta(x)=-\varphi'(x)
\end{eqnarray*}
for all $x\in(0,1)\cup(1,\infty)$.
Note that $\varphi\geq0$ is monotone decreasing. Then, by integration by parts, we find that
\begin{eqnarray*}
\begin{split}
F(t)
&=-\int\limits_0^\infty\varphi'(x)(\rho(0,t)-\rho(x,t))dx=\int\limits_0^\infty\varphi(x)(-\partial_x\rho(x,t))dx\\
&=\sum_{k\in\mathbb{Z}}\int\limits_{2^k}^{2^{k+1}}\varphi(x)(-\partial_x\rho(x,t))dx\leq\sum_{k\in\mathbb{Z}}\varphi(2^k)\int\limits_{2^k}^{2^{k+1}}(-\partial_x\rho(x,t))dx\\
&\leq\sum_{k\in\mathbb{Z}}\varphi(2^k)(\rho(2^k,t)-\rho(2^{k+1},t)).
\end{split}
\end{eqnarray*}
Furthermore, by the Cauchy-Schwarz inequality, we bound $F(t)$ as
\begin{eqnarray*}
\begin{split}
F(t)
&\leq\Big(\sum_{k\in\mathbb{Z}}\frac{(\varphi(2^k))^2}{\eta(2^{k+1})\log\Big(3\sqrt{\frac{2^{2k}+a^2}{9\cdot2^{2k}+a^2}}\Big)}\Big)^{\frac12}\\
&\ \ \ \ \ \ \
\times\Big(\sum_{k\in\mathbb{Z}}\eta(2^{k+1})\log\Big(3\sqrt{\frac{2^{2k}+a^2}{9\cdot2^{2k}+a^2}}\Big)(\rho(2^k,t)-\rho(2^{k+1},t))^2\Big)^{\frac12},
\end{split}
\end{eqnarray*}
which is equivalent to
\begin{eqnarray*}
&&\sum_{k\in\mathbb{Z}}\eta(2^{k+1})\log\Big(3\sqrt{\frac{2^{2k}+a^2}{9\cdot2^{2k}+a^2}}\Big)(\rho(2^{k+1},t)-\rho(2^k,t))^2\\
&&\ \ \ \
\geq\Big(\sum_{k\in\mathbb{Z}}\frac{(\varphi(2^k))^2}{\eta(2^{k+1})\log\Big(3\sqrt{\frac{2^{2k}+a^2}{9\cdot2^{2k}+a^2}}\Big)}\Big)^{-1}
(F(t))^2.
\end{eqnarray*}
This lower bound combining with \eqref{ODE-F-inequality} leads to
\begin{eqnarray*}
\frac{d}{dt}F(t)\geq\frac{g}{4\pi}\Big(\sum_{k\in\mathbb{Z}}\frac{(\varphi(2^k))^2}{\eta(2^{k+1})\log\Big(3\sqrt{\frac{2^{2k}+a^2}{9\cdot2^{2k}+a^2}}\Big)}\Big)^{-1}
(F(t))^2.
\end{eqnarray*}
In view of the choice of $\eta$ and $\varphi$, we have that
\begin{eqnarray*}
\begin{split}
&\sum_{k\in\mathbb{Z}}\frac{(\varphi(2^k))^2}{\eta(2^{k+1})\log\Big(3\sqrt{\frac{2^{2k}+a^2}{9\cdot2^{2k}+a^2}}\Big)}\\
&=\sum_{k<0}\frac{(\frac{1}{3+\alpha}+\frac{1}{1-\alpha}(1-2^{k(1-\alpha)}))^2}{2^{-\alpha(k+1)}\cdot\frac12\log\Big(\frac{9\cdot2^{2k}+9a^2}{9\cdot2^{2k}+a^2}\Big)}
+\sum_{k\geq0}\frac{\frac{1}{(3+\alpha)^2}\cdot2^{2k(-3-\alpha)}}{2^{-(k+1)(4+\alpha)}\cdot\frac12\log\Big(\frac{9\cdot2^{2k}+9a^2}{9\cdot2^{2k}+a^2}\Big)}\\
&\leq\frac{32(1-2^{-\alpha})^{-1}}{(1-\alpha)^2(3+\alpha)^2\log\frac{9+9a^2}{9+a^2}}
+\frac{2^{5+\alpha}}{(3+\alpha)^2}\sum_{k\geq0}\frac{2^{-(2+\alpha)k}}{\log\Big(1+\frac{8a^2}{9\cdot2^{2k}+a^2}\Big)}
\equiv c_{a,\alpha}<\infty,
\end{split}
\end{eqnarray*}
since the series $\displaystyle\sum_{k\geq0}\frac{2^{-(2+\alpha)k}}{\log\Big(1+\frac{8a^2}{9\cdot2^{2k}+a^2}\Big)}$ is convergent.
Therefore, we obtain
\begin{eqnarray*}
\frac{d}{dt}F(t)
\geq\frac{g}{4\pi c_{a,\alpha}}(F(t))^2,
\end{eqnarray*}
which implies a finite time blow-up for $F(t)$, since
\begin{eqnarray*}
F(0)=\int\limits_0^\infty \eta(x)(\rho_0(0)-\rho_0(x))dx>0
\end{eqnarray*}
for any $\rho_0$ which is not a constant. This contradicts \eqref{bounded-F} and thus completes the proof of Theorem \ref{singularity-formation-telescoping-sums}. \hfill$\square$
\section{The Kiselev-type inequality for the transform $H_a$}
In this section, we adapt the arguments in \cite{[Kiselev]} to generalize the ranges of $\sigma$ and $p$ in Proposition \ref{nonlinear-inequality-Kiselev} to $\sigma>-1$ and $p\geq1$.
Specifically, we will prove the Kiselev-type inequality for the transform $H_a$, whose analogous version for the Hilbert transform $H$ was first established by Kiselev in \cite{[Kiselev]}. We state the main result as
\begin{proposition}\label{Ha-inequality-Kiselev-p}
Let $f:\R\rightarrow\R$ be an even continuously differentiable function with $ f'\geq0$ on $(0,\infty)$ and $f(0)=0$. Then for every $a>0$, $p\geq1$ and $\sigma>0$, we have
\begin{eqnarray}\label{Ha-inequality-kiselev---}
-\int\limits^\infty_0\frac{H_a f(x) f'(x)(f(x))^{p-1}}{x^\sigma}dx
\geq
C_{p,\sigma}\int\limits^\infty_0\frac{a^2(f(x))^{p+1}}{x^{1+\sigma}(x^2+a^2)}dx.
\end{eqnarray}
The constant $C_{p,\sigma}$ may depend only on $p$ and $\sigma$. If the right hand side of \eqref{Ha-inequality-kiselev---} is infinite, the inequality is understood in the sense that the left hand side must also be infinite.
\end{proposition}

Before giving the proof of Proposition \ref{Ha-inequality-Kiselev-p}, we present a local in space upper bound for the transform $H_a$. For this purpose, we begin with a simple identity for the transform $H_a$ acting on even functions defined on $\R$.
\begin{lemma}\label{identity-Ha}
Let $f:\R\rightarrow\R$ be an even continuously differentiable function. Then for any $a>0$ and $x>0$, we have
\begin{eqnarray*}
H_a f(x)=\frac{1}{2\pi}\int\limits_{0}^{2x}f'(y)\log\frac{1+\frac{a^2}{x^2}}{1+\frac{a^2}{(x-y)^2}}dy+\frac{a^2}{\pi}\int\limits_x^{\infty}\frac{f(y-x)-f(x+y)}{y(y^2+a^2)}dy.
\end{eqnarray*}
\end{lemma}
\textbf{Proof.}
By \eqref{definition-Ha} and a simple change of variables, we derive that
\begin{eqnarray*}
\begin{split}
H_af(x)
&=\frac{a^2}{\pi}\Big(\int\limits_{-\infty}^{-x}+\int\limits_{-x}^{x}+\int\limits_x^{\infty}\Big)\frac{f(x-y)}{y(y^2+a^2)}dy\\
&=\frac{a^2}{\pi}\int\limits_x^{\infty}\frac{f(x-y)-f(x+y)}{y(y^2+a^2)}dy+\frac{1}{\pi}\int\limits_{0}^{2x}\frac{a^2f(y)dy}{(x-y)((x-y)^2+a^2)}.
\end{split}
\end{eqnarray*}
Furthermore, using \eqref{derivative-identity-3.2} along with integration by parts, we derive that
\begin{eqnarray*}
\begin{split}
\int\limits_{0}^{2x}\frac{a^2f(y)dy}{(x-y)((x-y)^2+a^2)}
&=-\int\limits_{0}^{2x}f(y)\frac{\partial}{\partial y}\log\Big(\frac{|y-x|\sqrt{x^2+a^2}}{x\sqrt{(y-x)^2+a^2}}\Big)dy\\
&=\int\limits_{0}^{2x}f'(y)\log\Big(\frac{|y-x|\sqrt{x^2+a^2}}{x\sqrt{(y-x)^2+a^2}}\Big)dy\\
&=\frac{1}{2}\int\limits_{0}^{2x}f'(y)\log\frac{1+\frac{a^2}{x^2}}{1+\frac{a^2}{(x-y)^2}}dy.
\end{split}
\end{eqnarray*}
Then we conclude the proof of Lemma \ref{identity-Ha}.  \hfill$\square$

As a consequence of Lemma \ref{identity-Ha}, we can derive a upper bound of the transform $H_a$ acting on the even and monotone increasing (on $(0,\infty)$) functions on $\R$.
\begin{coro}\label{upper-bound-Ha-Kiselev}
Fix $a>0$. Let $f:\R\rightarrow\R$ be an even continuously differentiable function with $ f'\geq0$ on $(0,\infty)$. Then for any $x>0$ and $q\in(1,2)$, we have
\begin{eqnarray*}
H_a f(x)\leq -\frac{a^2q(2-q)}{2\pi}\cdot\frac{f(qx)-f(q^{-1}x)}{x^2+a^2}.
\end{eqnarray*}
\end{coro}
\textbf{Proof.}
By virtue of the formula: for any $x\neq y$,
\begin{eqnarray*}
\frac{\partial}{\partial y}\log\frac{1+\frac{a^2}{x^2}}{1+\frac{a^2}{(x-y)^2}}=\frac{2a^2}{(y-x)((x-y)^2+a^2)},
\end{eqnarray*}
it is not difficult to check that the function
\begin{eqnarray*}
y\mapsto\log\frac{1+\frac{a^2}{x^2}}{1+\frac{a^2}{(x-y)^2}}
\end{eqnarray*}
is non-positive on $[0,x)\cup(x,2x]$, strictly decreasing on $[0,x)$ and strictly increasing on $(x,2x]$ with a singularity of negative infinity attained $y=x$. Therefore,
by Lemma \ref{identity-Ha} and the monotonicity of $f$ on $(0,\infty)$,
we can obtain that
\begin{eqnarray*}
\begin{split}
H_af(x)
&\leq\frac{1}{2\pi}\int\limits_{0}^{2x}f'(y)\log\frac{1+\frac{a^2}{x^2}}{1+\frac{a^2}{(x-y)^2}}dy
\leq\frac{1}{2\pi}\int\limits_{q^{-1}x}^{qx}f'(y)\log\frac{1+\frac{a^2}{x^2}}{1+\frac{a^2}{(x-y)^2}}dy\\
&\leq\frac{1}{2\pi}\max\Big\{\log\frac{1+\frac{a^2}{x^2}}{1+\frac{(aq)^2}{(q-1)^2x^2}},\log\frac{1+\frac{a^2}{x^2}}{1+\frac{a^2}{(q-1)^2x^2}}\Big\}\cdot\int\limits_{q^{-1}x}^{qx}f'(y)dy\\
&=\frac{1}{2\pi}\log\Big(1-\frac{q(2-q)a^2}{(q-1)^2x^2+a^2}\Big)(f(qx)-f(q^{-1}x))\\
&\leq\frac{1}{2\pi}\cdot\frac{-q(2-q)a^2}{(q-1)^2x^2+a^2}(f(qx)-f(q^{-1}x))\\
&\leq-\frac{a^2q(2-q)}{2\pi}\frac{f(qx)-f(q^{-1}x)}{x^2+a^2},
\end{split}
\end{eqnarray*}
which finishes the proof of Corollary \ref{upper-bound-Ha-Kiselev}. \hfill$\square$

With help of the above corollary, we are now ready to prove Proposition \ref{Ha-inequality-Kiselev-p}.

\textbf{Proof of Proposition \ref{Ha-inequality-Kiselev-p}.}
Fix a number $q$ such that $1<q<2$.
By Corollary \ref{upper-bound-Ha-Kiselev}, it is sufficient to prove the following inequality
\begin{eqnarray}\label{intermediate-inequality}
\begin{split}
\int\limits^\infty_0\frac{(f(qx)-f(q^{-1}x)) f'(x)(f(x))^{p-1}}{x^\sigma(x^2+a^2)}dx
&\geq C_{p,\sigma}\int\limits^\infty_0\frac{(f(x))^{p+1}}{x^{1+\sigma}(x^2+a^2)}dx
\end{split}
\end{eqnarray}
with $C_{p,\sigma}>0$ which may depend only on $p$ and $\sigma$. For this purpose, we first prove that for any given $L>0$, it holds that
\begin{eqnarray}\label{kiselev-L}
\begin{split}
\int\limits^L_0\frac{(f(qx)-f(q^{-1}x)) f'(x)(f(x))^{p-1}}{x^\sigma(x^2+a^2)}dx
\geq C'_{p,\sigma}\int\limits^L_0\frac{(f(x))^{p+1}}{x^{1+\sigma}(x^2+a^2)}dx,
\end{split}
\end{eqnarray}
where the constant $C'_{p,\sigma}$ depends only on $p,\sigma>0$, and is independent of $L>0$.

Indeed, split the integrals on both sides of \eqref{kiselev-L} into intervals $\Big[\frac{L}{q^{k+1}},\frac{L}{q^{k}}\Big],k\in\mathbb{N}$, and set $b_k\equiv f(q^{-k}L)$.
Notice that, for any $k=0,1,2,...$
\begin{eqnarray}\label{sk-lower-bound}
\begin{split}
l_k
&\equiv\int\limits^{q^{-k}L}_{q^{-k-1}L}\frac{(f(qx)-f(q^{-1}x)) f'(x)(f(x))^{p-1}}{x^\sigma(x^2+a^2)}dx\\
&\geq\frac{b_k-b_{k+1}}{q^{-2k}L^2+a^2}\int\limits^{q^{-k}L}_{q^{-k-1}L}\frac{ f'(x)(f(x))^{p-1}}{x^\sigma}dx.
\end{split}
\end{eqnarray}
Integration by parts in the last integral in \eqref{sk-lower-bound} directly yields that
\begin{eqnarray}\label{identity}
\int\limits^{q^{-k}L}_{q^{-k-1}L}\frac{ f'(x)(f(x))^{p-1}}{x^\sigma}dx
=\frac{q^{\sigma k}}{pL^\sigma}(b^p_k-b^p_{k+1})
+\frac{\sigma}{p}\int\limits^{q^{-k}L}_{q^{-k-1}L}\frac{ (f(x))^p-b^p_{k+1}}{x^{1+\sigma}}dx.
\end{eqnarray}
Also, for any $k=0,1,2,...$
\begin{eqnarray}\label{rk-upper-bound}
r_k\equiv\int\limits^{q^{-k}L}_{q^{-k-1}L}\frac{(f(x))^{p+1}dx}{x^{1+\sigma}(x^2+a^2)}
\leq\frac{(q-1)q^{\sigma( k+1)}}{L^{\sigma}(q^{-2k-2}L^2+a^2)}\cdot b^{p+1}_{k}.
\end{eqnarray}
In the spirit of the idea in \cite{[Kiselev]}, we introduce the so-called ``good" set $\mathcal{G}$ defined by
\begin{eqnarray*}
\mathcal{G}=\{k\in\mathbb{N}:b_{k}-b_{k+1}\geq cb_{k}\},
\end{eqnarray*}
where the fixed constant $c$ is chosen such that $(1-c)^{p+1}q^{\sigma }>1$.
We proceed to distinguish the proof into two cases.

\textbf{Case 1.} $\mathcal{G}$ is a finite set.
In this case, there exists an $K\in\mathbb{N}$ such that for all $k>K$, $b_{k+1}>(1-c)b_{k}$.
It follows that, for all $k>K$
\begin{eqnarray*}
\begin{split}
r_k
&\geq
\frac{(1-q^{-1})q^{\sigma k}b^{p+1}_{k+1}}{L^{\sigma}(q^{-2k}L^2+a^2)}
\geq\frac{(1-q^{-1})q^{\sigma k}}{L^{\sigma}(L^2+a^2)}\cdot((1-c)^{k-K}b_{K+1})^{p+1}\\
&=\frac{(1-q^{-1})b^{p+1}_{K+1}}{L^{\sigma}(L^2+a^2)(1-c)^{K(p+1)}}\cdot((1-c)^{p+1}q^{\sigma })^k\rightarrow\infty
\end{split}
\end{eqnarray*}
as $k\rightarrow\infty$ by our choice of $c$. This shows that
\begin{eqnarray*}
\int\limits^L_0\frac{(f(x))^{p+1}}{x^{1+\sigma}(x^2+a^2)}dx=\infty.
\end{eqnarray*}
By \eqref{sk-lower-bound} and \eqref{identity}, we continue to find that for all $k>{K}$,
\begin{eqnarray}\label{sk-infinity}
\begin{split}
l_k
&\geq\frac{b_k-b_{k+1}}{q^{-2k}L^2+a^2}\cdot \frac{q^{\sigma k}}{pL^{\sigma}}(b^p_k-b^p_{k+1})\\
&\geq\frac{(b_k-b_{k+1})^2}{L^2+a^2}\frac{q^{\sigma k}}{L^{\sigma}}b^{p-1}_{k+1}.
\end{split}
\end{eqnarray}
In order to continue to get the lower bound of $l_k$, we introduce $d_k$ such that $b_k-b_{k+1}=b_kd_k$, and $d_k<c$ if $k>K$.
Note that $\displaystyle\prod^\infty_{k=K+1}(1-d_k)=\displaystyle\lim_{k\rightarrow\infty}b_k=f(0)=0$.
Then we must have $\displaystyle\sum_{k}d_k=\infty$.
It follows from \eqref{sk-infinity} that
\begin{eqnarray*}
\begin{split}
l_k
&\geq \frac{q^{\sigma k}}{L^{\sigma}}\frac{d^2_kb^{p+1}_{k+1}}{L^2+a^2}\geq\frac{q^{\sigma k}}{L^{\sigma}}\frac{d^2_k}{L^2+a^2}((1-c)^{k-K}b_{K+1})^{p+1}\\
&=\frac{b^{p+1}_{K+1}}{L^{\sigma}(L^2+a^2)(1-c)^{K(p+1)}}\cdot ((1-c)^{p+1}q^{\sigma })^kd^2_k
\end{split}
\end{eqnarray*}
which along with $\displaystyle\sum_{k}d_k=\infty$ and the choice of $c$ implies that $\displaystyle\sum_{k}l_k=\infty$, that is,
\begin{eqnarray*}
\int\limits^L_0\frac{(f(qx)-f(q^{-1}x)) f'(x)(f(x))^{p-1}}{x^\sigma(x^2+a^2)}dx=\infty.
\end{eqnarray*}
In this case, the inequality \eqref{kiselev-L} holds for any constant $C'_{p,\sigma}>0$.

\textbf{Case 2.} $\mathcal{G}$ is an infinite set. In this case, we set
\begin{eqnarray*}
\mathcal{G}=\{k_j\in\mathbb{N}|j=1,2,3,\cdot\cdot\cdot\},
\end{eqnarray*}
where $k_j<k_{j+1}$ for any $j=1,2,3,\cdot\cdot\cdot$.
Clearly, the integral in the right-hand side in \eqref{kiselev-L} is
\begin{eqnarray*}
\sum^\infty_{k=0}r_k=\sum^\infty_{j=1}r_{k_j}+\sum_{m\in[0,k_1)}r_m+\sum^\infty_{j=2}\sum_{m\in(k_{j-1},k_j)}r_m.
\end{eqnarray*}
Note that when $k_{j}= k_{j-1}+1$, $(k_{j-1},k_j)\cap \mathcal{G}^c=\emptyset$, at this moment, we understand $r_m=0$ in the above summation. Also, when $k_1=0$, the second sum is $0$.

For $k\in\mathcal{G}$, by \eqref{sk-lower-bound}-\eqref{rk-upper-bound}, then we have that
\begin{eqnarray}\label{sk}
\begin{split}
l_k
&\geq\frac{b_k-b_{k+1}}{q^{-2k}L^2+a^2}\cdot \frac{q^{\sigma k}}{pL^{\sigma}}(b^p_k-b^p_{k+1})\geq\frac{cb_k}{q^{-2k}L^2+a^2}\cdot \frac{q^{\sigma k}}{pL^\sigma}(1-(1-c)^p)b^p_k\\
&\geq\frac{c^2r_k}{p(q-1)q^{\sigma+2}}\cdot\frac{q^{-2k}L^2+q^2a^2}{q^{-2k}L^2+a^2}
\geq\frac{c^2}{p(q-1)q^{\sigma+2}}\cdot r_k.
\end{split}
\end{eqnarray}
Also, for every $m\in(k_{j-1},k_{j})$, by \eqref{rk-upper-bound}, $(k_{j-1},k_{j})\cap\mathcal{G}=\emptyset$, $k_{j}\in\mathcal{G}$, and the second inequality in \eqref{sk}, we have that
\begin{eqnarray*}
\begin{split}
r_m
&\leq\frac{(q-1)q^{\sigma( m+1)}}{L^\sigma(q^{-2m-2}L^2+a^2)}b^{p+1}_m\\
&\leq\frac{(q-1)q^{\sigma( m+1)}}{L^\sigma(q^{-2m-2}L^2+a^2)}\Big(\frac{1}{1-c}\Big)^{(p+1)(k_j-m)}b^{p+1}_{k_j}\\
&\leq\frac{p(q-1)q^{\sigma+2}}{c^2}\Big(\frac{1}{(1-c)^{p+1}q^{\sigma}}\Big)^{k_j-m}\frac{q^{-2k_j}L^2+a^2}{q^{-2m}L^2+q^{2}a^2}l_{k_j}\\
&\leq\frac{p(q-1)q^{\sigma+2}}{c^2}\Big(\frac{1}{(1-c)^{p+1}q^{\sigma}}\Big)^{k_j-m}l_{k_j}.
\end{split}
\end{eqnarray*}
Thus,
\begin{eqnarray*}
\begin{split}
\sum^{k_j-1}_{m=k_{j-1}+1}r_{m}
&\leq\frac{p(q-1)q^{\sigma+2}}{c^2}l_{k_j}\sum^{k_j-1}_{m=k_{j-1}+1}\Big(\frac{1}{(1-c)^{p+1}q^{\sigma}}\Big)^{k_j-m}\\
&=\frac{p(q-1)q^{\sigma+2}}{c^2}l_{k_j}\sum^{k_j-k_{j-1}-1}_{m=1}\Big(\frac{1}{(1-c)^{p+1}q^{\sigma}}\Big)^{m}\\
&\leq \frac{p(q-1)q^{2\sigma+2}(1-c)^{p+1}}{c^2((1-c)^{p+1}q^{\sigma}-1)}l_{k_j}.
\end{split}
\end{eqnarray*}
Similarly, we can derive that, when $k_1>0$,
\begin{eqnarray*}
\sum_{m\in[0,k_1)}r_m\leq\frac{p(q-1)q^{2\sigma+2}(1-c)^{p+1}}{c^2((1-c)^{p+1}q^{\sigma}-1)}l_{k_1}.
\end{eqnarray*}
Finally, we have that
\begin{eqnarray*}
\begin{split}
\sum^\infty_{k=0}r_k
&\leq\frac{p(q-1)q^{\sigma+2}}{c^2}\sum^\infty_{j=1}l_{k_j}
+\frac{p(q-1)q^{2\sigma+2}(1-c)^{p+1}}{c^2((1-c)^{p+1}q^{\sigma}-1)}\sum^\infty_{j=1}l_{k_j}\\
&\leq\frac{p(q-1)q^{\sigma+2}}{c^2}\Big(1+\frac{q^{\sigma}(1-c)^{p+1}}{(1-c)^{p+1}q^{\sigma}-1}\Big)\sum^\infty_{k=0}l_k.
\end{split}
\end{eqnarray*}
This shows that \eqref{kiselev-L} holds true with the constant
\begin{eqnarray*}
C'_{p,\sigma}=\frac{c^2((1-c)^{p+1}q^{\sigma}-1)}{p(q-1)q^{\sigma+2}(2q^{\sigma}(1-c)^{p+1}-1)},
\end{eqnarray*}
which is independent of the length $L>0$.

Combining two cases we know the inequality \eqref{kiselev-L} is valid. Note that the constant $C'_{p,\sigma}$ in \eqref{kiselev-L} is independent of the parameter $L>0$. Letting $L\rightarrow\infty$ and using the Lebesgue Montone Convergence Theorem then yields the inequality \eqref{intermediate-inequality}. The proof of Proposition \ref{Ha-inequality-Kiselev-p} is finished.   \hfill$\square$

\appendix
\section{Proof of the local wellposedness}
In this appendix, we now sketch the proof of the local well-posedness  and the Beale-Kato-Majda type criterion for \eqref{KS}.
Before that, we recall that the transform $H_a$ satisfies the assumptions of Calder\'{o}n-Zygmund theory, with being a bounded linear operator on $L^p$ for any $p\in(1,\infty)$. In addition, the transform $H_a$ maps $L^\infty(\mathbb{R})$ to $BMO(\mathbb{R})$. We refer the readers to \cite{[Stein1967]} and \cite{[Stein1970]} for the proof of these properties for the singular integral operator $H_a$.

\textbf{Proof of Proposition \ref{local-well-posedness}.}
We first derive the local energy estimate.
Standard $L^2$ estimate gives that
\begin{eqnarray}\label{L2}
\begin{split}
\frac{d}{dt}\|\rho(t)\|_{L^2(\R)}
&\leq\frac12g\|\rho(t)\|_{L^2(\R)}\|H_a\partial_x\rho(t)\|_{L^\infty(\R)},
\end{split}
\end{eqnarray}
Note that, by the continuous embedding $H^\lambda(\mathbb{R})\hookrightarrow L^\infty(\mathbb{R})$ for $\lambda>\frac12$,
\begin{eqnarray*}
\begin{split}
|(H_a-H)\partial_x\rho(x,t)|
&\leq\frac{1}{\pi}\int\limits_{\R}\frac{|x-y|}{(x-y)^2+a^2}|\partial_y\rho(y,t)|dy\\
&\leq\frac{1}{\pi}\|\partial_x\rho(t)\|_{L^\infty}\int\limits_{|z|\leq1}\frac{|z|dz}{z^2+a^2}
+\frac{1}{\pi}\|\partial_x\rho(t)\|_{L^2}\Big(\int\limits_{|z|\geq1}\frac{z^2dz}{(z^2+a^2)^2}\Big)^{\frac12}\\
&\leq C_{a,s}\|\rho(t)\|_{H^s}.
\end{split}
\end{eqnarray*}
Then
\begin{eqnarray}\label{velocity-Lipschitz}
\begin{split}
\|H_a\partial_x\rho(t)\|_{L^\infty(\R)}
&\leq\|\Lambda\rho(t)\|_{L^\infty(\R)}+\|(H_a-H)\partial_x\rho(t)\|_{L^\infty(\R)}\\
&\leq C_{a,s}\|\rho(t)\|_{H^s}.
\end{split}
\end{eqnarray}
It follows that
\begin{eqnarray}\label{L2-local}
\begin{split}
\frac{d}{dt}\|\rho(t)\|_{L^2(\R)}
&\leq gC_{a,s}\|\rho(t)\|_{H^s}\|\rho(t)\|_{L^2(\R)}.
\end{split}
\end{eqnarray}
Next, applying $\Lambda^s:=(-\Delta)^{\frac s2}$ in \eqref{KS}, multiplying by $\Lambda^s\rho$ and integrating by parts, we can arrive at
\begin{eqnarray*}
\begin{split}
&\frac{1}{2}\frac{d}{dt}\|\Lambda^s\rho(t)\|^2_{L^2}
=-g\int\limits_{\R}\Lambda^s\rho\Big(\Lambda^s(H_a\rho\partial_x\rho)-H_a\rho\partial_x\Lambda^s\rho\Big) dx+\frac g2\int\limits_{\R}\partial_xH_a\rho(\Lambda^s\rho)^2 dx\\
&\leq g\|\Lambda^s\rho\|_{L^2}\|\Lambda^s(H_a\rho\partial_x\rho)-H_a\rho\partial_x\Lambda^s\rho\|_{L^2}+\frac g2\|H_a\partial_x\rho\|_{L^\infty}\|\Lambda^s\rho\|^2_{L^2}.
\end{split}
\end{eqnarray*}
By the classical commutator estimate of the operator $\Lambda^s$ (see, e.g., \cite{[Ju]}), we can estimate
\begin{eqnarray*}
\begin{split}
\|\Lambda^s(H_a\rho\partial_x\rho)-H_a\rho\partial_x\Lambda^s\rho\|_{L^2}
&\leq C_s\Big(\|\Lambda^sH_a\rho\|_{L^2}\|\partial_x\rho\|_{L^\infty}
+\|\Lambda^{s-1}\partial_x\rho\|_{L^2}\|\partial_x H_a\rho\|_{L^\infty}\Big)\\
&\leq C_{a,s}\Big(\|\partial_x\rho\|_{L^\infty}
+\|\partial_x H_a\rho\|_{L^\infty}\Big)\|\Lambda^s\rho\|_{L^2}.
\end{split}
\end{eqnarray*}
Therefore, we have that
\begin{eqnarray}\label{High-order-derivative-estimate}
\frac{d}{dt}\|\Lambda^s\rho(t)\|_{L^2(\R)}
\leq gC_{a,s}\Big(\|\partial_x\rho\|_{L^\infty}
+\|\partial_x H_a\rho\|_{L^\infty}\Big)\|\Lambda^s\rho\|_{L^2}..
\end{eqnarray}
Finally, inserting \eqref{velocity-Lipschitz} into \eqref{High-order-derivative-estimate}, and adding with \eqref{L2-local}, we conclude that
\begin{eqnarray*}
\frac{d}{dt}\|\rho(t)\|_{H^s(\R)}
\leq gC_{a,s}\|\rho(t)\|^2_{H^s(\R)}.
\end{eqnarray*}
In order to make the above formal energy estimates rigorous, one needs to work with the regularized system
\begin{eqnarray*}
\partial_t\rho+gJ_\epsilon(H_aJ_\epsilon\rho\partial_xJ_\epsilon\rho)=0,
\end{eqnarray*}
where $J_\epsilon$ is the usual mollifier. We refer the interested reader to \cite{[Chae],[Maida-Bertozzi]} for more details.

We proceed to prove the uniqueness part. Suppose $\rho_1,\rho_2\in C([0,T); H^s(\R))$ are solutions to \eqref{KS} with the same initial data $\rho_0\in H^s$.
By \eqref{KS}, we have
\begin{eqnarray*}
\partial_t(\rho_1-\rho_2)=-gH_a\rho_1\partial_x(\rho_1-\rho_2)-gH_a(\rho_1-\rho_2)\partial_x\rho_2.
\end{eqnarray*}
It follows that
\begin{eqnarray*}
\begin{split}
\frac{d}{dt}\|\rho_1-\rho_2\|_{L^2}
&\leq\frac g2\|\partial_x H_a\rho_1\|_{L^\infty}\|\rho_1-\rho_2\|_{L^2}
+g\|H_a(\rho_1-\rho_2)\|_{L^2}\|\partial_x\rho_2\|_{L^\infty}\\
&\leq gC_{a,s}(\|\rho_1\|_{H^s}+\|\rho_2\|_{H^s})\|\rho_1-\rho_2\|_{L^2}.
\end{split}
\end{eqnarray*}
Gr\"{o}nwall's inequality along with $(\rho_1-\rho_2)(x,0)=0$ finishes the proof of the uniqueness of solutions.

Finally, we prove the Beale-Kato-Majda type criterion for \eqref{KS}.
Indeed, by \eqref{L2} and \eqref{High-order-derivative-estimate}, we have
\begin{eqnarray}\label{3.10}
\frac{d}{dt}\|\rho(t)\|_{H^s}\leq g C_{a,s}(\|\partial_x\rho\|_{L^\infty}+\|\partial_xH_a\rho\|_{L^\infty})\|\rho(t)\|_{H^s}.
\end{eqnarray}
With the help of the limiting Sobolev inequality (see, e.g., Theorem 1 in \cite{[Kozono-Taniuchi]}): for any $f\in H^\lambda(\R)$ with $\lambda>\frac12$,
\begin{eqnarray*}
\|f\|_{L^\infty}\leq C_\lambda(1+\|f\|_{BMO})(1+\log^+\|f\|_{H^\lambda}),
\end{eqnarray*}
we have that, for $s>\frac32$,
\begin{eqnarray*}
\begin{split}
\|\partial_xH_a\rho\|_{L^\infty}
&\leq C_s(1+\|H_a\partial_x\rho\|_{BMO})(1+\log^+\|H_a\partial_x\rho\|_{H^{s-1}})\\
&\leq C_{a,s}(1+\|\partial_x\rho\|_{L^\infty})\log(e+\|\rho\|_{H^{s}}).
\end{split}
\end{eqnarray*}
Substituting this logarithmic-type estimate into \eqref{3.10}, we obtain that
\begin{eqnarray*}
\frac{d}{dt}\|\rho(t)\|_{H^s}\leq C_{a,s}(1+\|\partial_x\rho\|_{L^\infty})\log(e+\|\rho\|_{H^{s}})\|\rho(t)\|_{H^s}.
\end{eqnarray*}
This completes the proof of the blow-up criterion by using Gr\"{o}nwall's inequality.
The proof of Proposition \ref{local-well-posedness} is now finished.     \hfill$\square$

{\bf Acknowledgements.}
W. Zhang was supported by the Science and Technology Research Project of Department of Education of Jiangxi Province, China (No. GJJ2200363).

\end{document}